\documentclass{amsart}

\usepackage{amsmath,amssymb,amsthm,amscd}

\newcommand{\ot}{\otimes}
\newcommand{\fand}{\quad\text{and}\quad}
\newcommand{\ts}{\,}
\newcommand{\tss}{\hspace{1pt}}
\newcommand{\st}{\widetilde{\si}}
\newcommand{\Tb}{\widehat{T}}
\newcommand{\Tt}{\widetilde{T}}

\newcommand{\U}{{\rm U}}
\newcommand{\Ar}{{\rm A}}
\newcommand{\Sr}{{\rm S}}
\newcommand{\Y}{{\rm Y}}
\newcommand{\End}{{\rm{End}\ts}}
\newcommand{\gr}{{\rm gr}\ts}
\newcommand{\grpr}{{\rm gr}^{\tss\prime}\ts}
\newcommand{\degpr}{{\rm deg}^{\tss\prime}\tss}
\newcommand{\DY}{{\rm DY}}
\newcommand{\id}{{\rm id}}

\newcommand{\CC}{\mathbb{C}}
\newcommand{\ZZ}{\mathbb{Z}}

\newcommand{\Rc}{{\mathcal R}}
\newcommand{\Sc}{{\mathcal S}}

\newcommand{\Sym}{\mathfrak S}
\newcommand{\gl}{\mathfrak{gl}}
\newcommand{\g}{\mathfrak{g}}
\newcommand{\agot}{\mathfrak{a}}

\newcommand{\sll}{\mathfrak{sl}}

\newcommand{\al}{\alpha}
\newcommand{\be}{\beta}

\newcommand{\de}{\delta}
\newcommand{\De}{\Delta}
\newcommand{\ve}{\varepsilon}
\newcommand{\si}{\sigma}
\newcommand{\om}{\theta}

\newtheorem{thm}{Theorem}[section]
\newtheorem{lemma}[thm]{Lemma}
\newtheorem{prop}[thm]{Proposition}
\newtheorem{cor}[thm]{Corollary}

\newcommand{\bth}{\begin{thm}}
\renewcommand{\eth}{\end{thm}}
\newcommand{\bpr}{\begin{prop}}
\newcommand{\epr}{\end{prop}}
\newcommand{\ble}{\begin{lemma}}
\newcommand{\ele}{\end{lemma}}
\newcommand{\bco}{\begin{cor}}
\newcommand{\eco}{\end{cor}}

\newcommand{\bpf}{\begin{proof}}
\newcommand{\epf}{\end{proof}}

\def\beql#1{\begin{equation}\label{#1}}

\newcommand{\bal}{\begin{aligned}}
\newcommand{\eal}{\end{aligned}}
\newcommand{\beq}{\begin{equation}}
\newcommand{\eeq}{\end{equation}}
\newcommand{\ben}{\begin{equation*}}
\newcommand{\een}{\end{equation*}}

\numberwithin{equation}{section}

\begin{document}

\title{Double Yangian and the universal {\textit{R\tss}}-matrix}

\author{Maxim Nazarov}




\begin{abstract}
We describe the double Yangian of the general linear 
Lie algebra $\gl_N$ by following a general scheme of Drinfeld.
This description is based on the construction of the universal 
$R$-matrix for the Yangian.
To make exposition self contained, we include the proofs of
all necessary facts about the Yangian itself. In particular, we 
describe the centre
of the Yangian by using its Hopf algebra structure, 
and provide a 
proof of the analogue of the
Poincar\'e--Birkhoff--Witt theorem for the Yangian 
based on its representation theory. 
This proof extends to the double Yangian,
thus giving a description of its underlying vector space.
\end{abstract}


\maketitle

\tableofcontents

\thispagestyle{empty}



\section*{Introduction}

The main subject of this article is a Hopf algebra that 
appeared in the framework of quantum inverse scattering method 
introduced by L.\,D.\,Faddeev, E. K. Sklyanin and their collaborators, 
see for instance 
\cite{f:im,ks,s,stf,s:qv,tf}. 
This algebra then 
became a part of a family
of examples in the theory of quantum groups created by V.\,G.\,Drinfeld
\cite{d:ha,d:nr,d:qg}. He gave to this family the name
\textit{Yangians\/} in honour of C.\,N.\,Yang, 
the author of a seminal work \cite{y}. The Yangian that we
consider here 
corresponds to the general linear Lie algebra $\gl_N\,$.
It is a canonical deformation of the universal enveloping
algebra of the polynomial current Lie algebra $\gl_N[z]\,$.
 
The general notion of a quantum double was also introduced in \cite{d:qg}.
However the Yangians were not discussed there in the context of this notion. 
Here we define the double Yangian of the Lie algebra 
$\gl_N$ similarly to \cite{i:br}. Yet many details and proofs 
are also missing in the latter work.
In the present article we fill these gaps.

We denote by $\Y(\gl_N)$ the Yangian of 
$\gl_N\,$, and by $\DY(\gl_N)$ its quantum double.
There are several equivalent definitions of the Hopf
algebra $\Y(\gl_N)$ available \cite{m}. 
In this article we use the definition that 
appeared first, see for instance \cite{ks:sy,ks:qs,ft:ss}.
Details of this definition are given in our Sections
\ref{sec:mot},\ref{sec:matrix} and \ref{sec:hopf} 
by closely following~\cite{mno:yc}.
Sections \ref{sec:auto},\ref{sec:filtr} and \ref{sec:provec}
describe basic properties of the Yangian $\Y(\gl_N)$
that we will use. 

We will also use an analogue of the classical
Poincar\'e--Birkhoff--Witt theorem \cite{d:ae}
for the algebra $\Y(\gl_N)\,$. The first proof
of this analogue was given by V. G. Drinfeld but not published.
Other proofs were given later in \cite{l:pb,o:ri}.
In Section \ref{sec:pbw} we give yet another proof of this analogue 
by using the representation theory of current~Lie algebras. 
The fact from the theory that we use is established in 
Section~\ref{sec:eval}. It~is~this proof that will be extended 
to the double Yangian $\DY(\gl_N)$ in the present article.
This method was used in \cite{n:yq}
to prove analogues of the Poincar\'e--Birkhoff--Witt theorem 
for the Yangian of the queer Lie superalgebra $\mathfrak{q}_N$ and 
its quantum double. For the algebra dual to the coalgebra
$\Y(\gl_N)$ the same method was used in \cite{ek3}.

The structure of a Hopf algebra includes a canonical anti-automorphism
relative to both multiplication and comultiplication, called the antipodal 
map. In general this map is not involutive. In Section \ref{sec:auto} 
we also compute the square of this~map for the Yangian $\Y(\gl_N)\,$,
by following \cite{n:qb} where the Yangian of the general linear Lie 
superalgebra $\mathfrak{gl}_{\tss M|N}$ was considered.
This yields a family of central 
elements of the algebra $\Y(\gl_N)\,$, see also \cite{d:acha}.
In Section \ref{sec:cen} 
we prove that these elements generate 
the whole centre. 
Our proof uses another 
general fact from the theory
of current Lie algebras, 
which we establish in the beginning 
of the section. 
The idea of reducing the proof to that fact
belongs to V. G. Drinfeld, as acknowledged in \cite{mno:yc}.

In Section \ref{sec:dualyang}
we introduce the bialgebra $\Y^\ast(\gl_N)$ dual to $\Y(\gl_N)\,$.
First we define it in terms generators and relations
similarly to $\Y(\gl_N)\,$. However $\Y^\ast(\gl_N)$ is not a Hopf algebra.
The antipodal map is defined only on a certain completion
$\Y^\circ(\gl_N)$ of $\Y^\ast(\gl_N)$ described at the end of that section. 
In Section \ref{sec:pairing} we define
a bialgebra pairing 
between $\Y(\gl_N)$ and $\Y^\ast(\gl_N)\,$.
This definition goes back to \cite{rtf:ql} where the
quantized universal enveloping algebras 
of simple Lie algebras were considered. 
In Section \ref{sec:nondeg} we prove
that this pairing is non-degenerate.
Details of this proof first appeared in \cite{n:yq} 
where instead of $\gl_N\,$,
the Lie superalgebra $\mathfrak{q}_N$
was considered.

In Section \ref{sec:unir} we define the universal $R$-matrix for
$\Y(\gl_N)$. This is the canonical element of a suitable
completion of the tensor product $\Y^\ast(\gl_N)\ot\Y(\gl_N)\,$,
which corresponds to the bialgebra pairing. There
we also describe the basic properties of this element
relative to the Hopf algebra structures on both 
$\Y(\gl_N)$ and $\Y^\circ(\gl_N)\,$.

\enlargethispage{4pt}

In Section \ref{sec:doubleyang} we define the double Yangian $\DY(\gl_N)$
as a bialgebra generated by $\Y(\gl_N)$ and $\Y^\ast(\gl_N)\,$.
Following 
\cite{d:qg,rtf:ql} the cross relations between the elements of $\Y(\gl_N)$ 
and $\Y^\ast(\gl_N)$ are introduced by means of the universal $R$-matrix.
Then~we provide a more explicit description of the algebra 
$\DY(\gl_N)\,$. Using this description one can define a central extension
of $\DY(\gl_N)\,$, see for instance \cite{ek45,jkmy}.

Finally, in Section \ref{sec:dpbw} we introduce a filtration on the algebra
$\DY(\gl_N)$ and show that the corresponding graded algebra
is isomorphic to the universal enveloping
algebra of the current Lie algebra $\gl_N[z,z^{-1}]\,$.
This implies our
analogue of the Poincar\'e--Birkhoff--Witt theorem for
$\DY(\gl_N)\,$. This also implies that 
the defining homomorphisms of the algebras
$\Y(\gl_N)$ and $\Y^\ast(\gl_N)$ to $\DY(\gl_N)$ are embeddings.

The purpose of the present article is to provide the basic
facts about the double Yangian $\DY(\gl_N)$ with their detailed
proofs. We do not not aim to review all works which involve this remarkable
object. Still let us mention here the pioneering works
\cite{bl:qd,ls:iq,sm} where the double Yangian of the special linear 
Lie algebra $\sll_2$ was studied. 
The double Yangians of all simple Lie 
algebras were studied in \cite{k,kt} by using the definition of the 
underlying Yangians from \cite{d:nr}.
This approach to double Yangians is different from ours.
Recently some of the results on $\DY(\gl_N)$ presented here
have been 
extended to the double Yangians 
of the other classical Lie algebras \cite{jly}.


\section{Definition of the Yangian}
\label{sec:mot}

The {\em Yangian\/} of the general linear Lie algebra 
$\gl_N$
is a unital associative algebra $\Y(\gl_N)$ 
over the complex field $\CC$
with countably many generators 
$$
T_{ij}^{(1)},\ T_{ij}^{(2)},\ldots
\quad\text{where}\quad
i,j=1,\dots,N\,.
$$ 
The defining relations of the algebra $\Y(\gl_N)$ are
\beql{defyang}
[\,T^{(r+1)}_{ij},T^{(s)}_{kl}\,]-[\,T^{(r)}_{ij}, T^{(s+1)}_{kl}\,]=
T^{(r)}_{kj} T^{(s)}_{il}-T^{(s)}_{kj} T^{(r)}_{il}
\eeq
where $r,s=0,1,\dots$ and  $T^{(0)}_{ij}=\delta_{ij}\,$.
By introducing the formal generating series
\beql{tiju}
T_{ij}(u) = \delta_{ij} + T^{(1)}_{ij} u^{-1} + T^{(2)}_{ij}u^{-2} +
\ldots\in\Y(\gl_N)[[u^{-1}]]
\eeq
we can write \eqref{defyang} in the form
\beql{defrel}
(u-v)\,[\,T_{ij}(u),T_{kl}(v)\,]=
T_{kj}(u)\ts T_{il}(v)-T_{kj}(v)\ts T_{il}(u)\,.
\eeq
Here the indeterminates $u$ and $v$ are considered to be commuting
with each other and with the elements of the Yangian.
The following is an equivalent form of 
\eqref{defyang}.

\bpr\label{prop:defequiv}
The system of relations \eqref{defyang} is equivalent to the system
\beql{defequiv}
[\,T^{(r)}_{ij}, T^{(s)}_{kl}\,]=\sum_{a=1}^{\min(r,s)}
\Big(T^{(a-1)}_{kj} T^{(r+s-a)}_{il}-T^{(r+s-a)}_{kj} T^{(a-1)}_{il}\Big)\,.
\eeq
\epr

\bpf
Observe that
the multiplication of both sides of \eqref{defrel} by the formal
series $\sum_{p=0}^{\infty} u^{-p-1}v^p$ yields an equivalent relation
\ben
[\,T_{ij}(u),T_{kl}(v)\,]=
\Big(T_{kj}(u)T_{il}(v)-T_{kj}(v)T_{il}(u)\Big)\sum_{p=0}^{\infty}u^{-p-1}v^p.
\een
Taking the coefficients of $u^{-r}v^{-s}$ on both sides gives
\ben
[\,T^{(r)}_{ij}, T^{(s)}_{kl}\,] =\sum_{a=1}^{r}
\Big(T^{(a-1)}_{kj} T^{(r+s-a)}_{il}-T^{(r+s-a)}_{kj} T^{(a-1)}_{il}\Big)\,.
\een
This agrees with \eqref{defequiv} in the case $r\leqslant s$. Finally, if
$r>s$ observe that
\ben
\sum_{a=s+1}^{r}
\Big(T^{(a-1)}_{kj} T^{(r+s-a)}_{il}-T^{(r+s-a)}_{kj} T^{(a-1)}_{il}\Big)=0\,.
\qedhere
\een
\epf

We shall be often using formal series to define or describe
maps between various algebras. If $A(u)$ and $B(u)$ are formal
series in $u$ with coefficients in certain algebras then
assignments of the type $A(u)\mapsto B(u)$
are understood in the sense that every coefficient of $A(u)$
is mapped to the corresponding coefficient of $B(u)$.

Many applications of $\Y(\gl_N)$
are based on the following observation.
Let $E_{ij}$ be the standard generators of the Lie algebra $\gl_N$ so that
\beq
\label{deflie}
[\,E_{ij},E_{kl}\,]=\de_{jk}\,E_{il}-\de_{li}\,E_{kj}\,.
\eeq

\bpr\label{prop:eval}
The assignment
\beql{eval}
T_{ij}(u)\mapsto \delta_{ij}+E_{ij}u^{-1}
\eeq
defines a surjective homomorphism $\Y(\gl_N)\to\U(\gl_N)\,$.
The assignment
\beql{embed}
E_{ij}\mapsto T_{ij}^{(1)}
\eeq
defines an embedding $\U(\gl_N)\to\Y(\gl_N)$.
\epr

\bpf
By the definition~\eqref{defrel} we need to verify the equality
\begin{gather*}
(u-v)\,[\,E_{ij}, E_{kl}\,]\ts u^{-1}v^{-1}=
\\
(\delta_{kj}+E_{kj}u^{-1})(\delta_{il}+E_{il}v^{-1})
-(\delta_{kj}+E_{kj}v^{-1})(\delta_{il}+E_{il}u^{-1})\,.
\end{gather*}
But this clearly holds by the commutation relations \eqref{deflie}
in $\gl_N$, which proves the first part of the proposition.
In order to prove the second part, put $r=s=1$ in \eqref{defequiv}.
This gives
\ben
[\,T_{ij}^{(1)}, T_{kl}^{(1)}\,]=
\de_{kj}T_{il}^{(1)}-\de_{il}T_{kj}^{(1)}.
\een
Thus \eqref{embed} is an algebra homomorphism. Its injectivity
follows from the observation that by applying \eqref{embed}
and then \eqref{eval}, we get the identity map on $\U(\gl_N)\,$.
\epf

The homomorphism \eqref{eval} is called
the \textit{evaluation homomorphism\/}. By its virtue
any representation of the Lie algebra $\gl_N$ can be
regarded as representation of the $\Y(\gl_N)\tss$.
Any irreducible representation of $\gl_N$ remains
irreducible over $\Y(\gl_N)$ due to surjectivity of this homomorphism.
We will also be using its composition
with the automorphism $E_{ij}\mapsto -E_{ji}$ of the algebra
$\U(\gl_N)\tss$. The composition maps
\beql{pin}
T_{ij}(u)\mapsto\delta_{ij}-E_{ji}u^{-1}\,.
\eeq
The reason for using it rather than \eqref{eval}
will be explained in Section~\ref{sec:provec}.


\section{Matrix form of the definition}
\label{sec:matrix}

Introduce the $N\times N$ matrix $T(u)$ whose $ij$-th entry is the series
$T_{ij}(u)\,$. One can regard $T(u)$ as an element
of the algebra $\End \CC^N\ot \Y(\gl_N)[[u^{-1}]]\,$. Then
\beql{tutens}
T(u)=\sum_{i,j=1}^N e_{ij}\ot T_{ij}(u)
\eeq
where $e_{ij}\in\End\CC^N$ are the standard matrix units. 
If $e_1,\dots,e_N$ are the standard basis vectors of
$\CC^N$, then  $T(u)\ts e_j$ is interpreted as the linear combination
\ben
T(u)\ts e_j=\sum_{i=1}^N e_i\ot T_{ij}(u)\in \CC^N\ot\Y(\gl_N)[[u^{-1}]]\,.
\een

For any positive integer $m$ we shall be using algebras of the form
\beql{multitp}
(\End\CC^N)^{\ot m}\ot\Y(\gl_N)\,.
\eeq
For any $a=1,\dots,m$ we denote
by $T_a(u)$
the matrix $T(u)$ which corresponds
to the $a$-th copy of the algebra $\End \CC^N$ in the tensor product
\eqref{multitp}. That is,
$T_a(u)$ is a formal power series in $u^{-1}$
with the coefficients from the algebra
\eqref{multitp},
$$
T_a(u)=\sum_{i,j=1}^N
1^{\ot(a-1)}\ot e_{ij}\ot 1^{\ot(m-a)}\ot T_{ij}(u)
$$
where 
$e_{ij}$ belongs to the $a$-th copy of $\End\CC^N$ and
$1$ is the identity matrix. If $C$ is an element 
of the tensor square $(\End\CC^N)^{\ot 2}$ 
then for $a,b=1,\dots,m$ with $a<b$ we will denote by $C_{ab}$
the image of $C$ under this embedding
$(\End\CC^N)^{\ot 2}\to(\End\CC^N)^{\ot m}:$
$$
e_{ij}\ot e_{kl}\mapsto
1^{\ot(a-1)}\ot e_{ij}\ot 1^{\ot(b-a-1)}\ot e_{kl}\ot 1^{\ot(m-b)}.
$$
Here the tensor factors 
$e_{ij}$ and $e_{kl}$ belong to
the $a$-th and $b$-th copies of $\End \CC^N$ respectively.
The element $C_{ab}$
can be identified with the element $C_{ab}\ot1$ of 
\eqref{multitp}.~If
\ben
t:\End\CC^N\to\End\CC^N:e_{ij}\mapsto e_{ji}
\een
is the matrix transposition, then
for any $a=1,\dots,m$ we shall denote by $t_a$ the corresponding
partial
transposition on the algebra \eqref{multitp}. It acts as $t$ on the
$a$-th copy of
$\End \CC^N$ and as the identity map on all the other tensor factors.

Consider now the permutation operator
\beql{pmatrix}
P=\sum_{i,j=1}^N e_{ij}\ot e_{ji}\in \End \CC^N\ot\End \CC^N.
\eeq
The rational function
\beql{rmatrix}
R(u)=1-P\,u^{-1}
\eeq
with values in $\End \CC^N\ot\End \CC^N$ is called
the {\it Yang\/} $R$-{\it matrix\/}.
Here and below we write $1$ instead of
$1\ot 1\tss$, for brevity.
We will be frequently using the identity
\ben
R(u)\tss R(-u)=1-u^{-2}.
\een

We will also work with the rational function
$$
R^{\tss t}(u)=1-{Q}\ts{u^{-1}}
$$
where
$$
Q=\sum_{i,j=1}^N e_{ij}\ot e_{ij}=P^{\,t_1}=P^{\,t_2}.
$$
We should write either
$R^{\tss t_1}(u)$ or $R^{\tss t_2}(u)$ instead of $R^{\tss t}(u)$
but we will \textit{not\/} do so.
Note that $Q$ is a one-dimensional operator
on $\CC^N\ot\CC^N$ such that $Q^2=N\ts Q\,$.
Hence
\beql{rtinv}
R^{\tss t}(u)^{-1}=1+Q\ts(u-N)^{-1}.
\eeq

\begin{prop}\label{prop:ybe}
In the algebra $(\End \CC^N)^{\ot3}(u,v)$
we have the identity
\beql{ybe}
R_{12}(u)\ts R_{13}(u+v)\ts R_{23}(v)=
R_{23}(v)\ts R_{13}(u+v)\ts R_{12}(u).
\eeq
\end{prop}

\bpf
Multiplying both sides of the relation
\eqref{ybe} by $uv(u+v)$ we come to verify
\beql{ybepol}
(u+P_{12})(u+v+P_{13})(v+P_{23})=(v+P_{23})(u+v+P_{13})(u+P_{12}).
\eeq
Each operator $P_{ij}$ is the image of the
corresponding transposition $(ij)\in\Sym_3$
under the natural action of the symmetric group
$\Sym_3$ on $(\CC^N)^{\ot 3}$
by permutations of the tensor factors. So \eqref{ybepol}
follows from the relations in the group algebra $\CC[\Sym_3]$.
\epf

The relation \eqref{ybe} is known as the
{\it Yang--Baxter equation}. The Yang $R$-matrix
is its simplest nontrivial solution.
Below we regard
$T_1(u)$ and $T_2(v)$ as formal power series
with the coefficients from the algebra
\eqref{multitp} where $m=2\,$. 
We also identify $R(u-v)$ with the rational function $R(u-v)\ot1$
taking values in this algebra.

\bpr\label{prop:ternary}
The defining relations of the algebra $\Y(\gl_N)$
can be written as
\beql{ternary}
R(u-v)\ts T_1(u)\ts T_2(v)=T_2(v)\ts T_1(u)\ts R(u-v).
\eeq
\epr

\bpf
Let us apply both sides of 
\eqref{ternary} to
an any basis vector $e_j\ot e_l\in\CC^N\ot\CC^N$
as explained in the beginning of this section.
For the left hand side we get
\ben
\sum_{i,k} T_{ij}(u)\ts T_{kl}(v)\ot e_i\ot e_k-
\frac{1}{u-v}\sum_{i,k}T_{ij}(u)\ts T_{kl}(v)\ot e_k\ot e_i\,,
\een
while the right hand side gives
\ben
\sum_{i,k} T_{kl}(v)\ts T_{ij}(u)\ot e_i\ot e_k-
\frac{1}{u-v}\sum_{i,k}T_{kj}(v)\ts T_{il}(u)\ot e_i\ot e_k\,.
\een
Multiplying by $u-v$ and
equating the coefficients of $e_i\ot e_k$ we recover
\eqref{defrel}.
\epf


\section{Automorphisms and anti-automorphisms}
\label{sec:auto}

In this section, we will use the $N\times N$ matrix $T(u)$
to define several distinguished automorphisms and anti-automorphisms
of the associative unital algebra $\Y(\gl_N)\tss$.
For each of them, we will
describe the $N\times N$ matrix whose $ij$-entry
is the formal power series in $u^{-1}$ with the coefficients
being the images of the corresponding coefficients of
the series $T_{ij}(u)$. For example, the assignment
\eqref{sign} below means that
for all indices $r=1,2,\dots$ and $i,j=1,\dots,N$
$$
T_{ij}^{(r)}\mapsto(-1)^r\,T_{ij}^{(r)}\,.
$$

\bpr\label{prop:auto}
For any $c\in\CC$ an automorphism of $\Y(\gl_N)$ can be defined by 
\beql{shift}
T(u)\mapsto T(u-c)\,.
\eeq
\epr

\bpf
The image of $T(u)$ relative to \eqref{shift} clearly
satisfies the defining relation~\eqref{ternary}.
Further, the mapping \eqref{shift} is obviously invertible
which completes the proof.
\epf

We may regard the element $T(u)$ defined by \eqref{tutens}
as a formal power series in $u^{-1}$
whose coefficients are matrices with the entries from
the algebra $\Y(\gl_N)$. Since the leading term of this series
is the identity matrix, the element $T(u)$ is invertible.
We denote by $T^{-1}(u)$ the inverse element.
Further, denote by $T^{\tss t}(u)$ the transposed matrix for $T(u)$. Then
\ben
T^{\tss t}(u)=\sum_{i,j=1}^N e_{ij}\ot T_{ji}(u)\,.
\een

\bpr
\label{prop:anti}
Each of the assignments
\begin{align}
\label{sign}
T(u)&\mapsto T(-u)\,,
\\
\label{transp}
T(u)&\mapsto T^{\tss t}(u)\,,
\\
\label{inverse}
\Sr:T(u)&\mapsto T^{-1}(u)
\end{align}
defines an anti-automorphism of\/ $\Y(\gl_N)\,$.
\epr

\bpf
The images $T_{ij}^{\tss\prime}(u)$ of the series $T_{ij}(u)$
under any anti-automorphism of the algebra $\Y(\gl_N)$ must satisfy
the relations
\eqref{defrel} with the opposite multiplication:
\ben
(u-v)\ts [\,T_{ij}^{\tss\prime}(u),T_{kl}^{\tss\prime}(v)\,]=
T_{il}^{\tss\prime}(u)\,T_{kj}^{\tss\prime}(v)-
T_{il}^{\tss\prime}(v)\,T_{kj}^{\tss\prime}(u).
\een
Exactly as in the proof of Proposition~\ref{prop:ternary}, one can show
that these relations can be equivalently written
in the following matrix form
$$
R(u-v)\ts T_2^{\tss\prime}(v)\ts T_1^{\tss\prime}(u) =
T_1^{\tss\prime}(u)\ts T_2^{\tss\prime}(v)\ts  R(u-v)
$$
where $T^{\tss\prime}(u)$ is the $N\times N$ matrix whose $ij$-th entry is
$T_{ij}^{\tss\prime}(u)$. 
But the relation
$$
R(u-v)\ts T_2(-v)\ts T_1(-u) =
T_1(-u)\ts T_2(-v)\ts  R(u-v)
$$
follows from \eqref{ternary} if we conjugate both sides by $P$ and
replace $(u,v)$ by $(-v,-u)\,$. This shows that \eqref{sign}
defines an anti-homomorphism. Furthermore, the application of the
partial transposition $t_1$ to both sides of the
relation \eqref{ternary} yields
\beql{trter}
T^{\tss t}_1(u)\ts R^{\tss t}(u-v)\ts T_2(v)  =
T_2(v)\ts R^{\tss t}(u-v)\ts T^{\tss t}_1(u).
\eeq
Since $R(u-v)$ is fixed by the composition of $t_1$ with $t_2\,$,
applying $t_2$ to \eqref{trter} yields
\ben
T^{\tss t}_1(u)\ts  T^{\tss t}_2(v)\ts R(u-v)  =
R(u-v)\ts T^{\tss t}_2(v)\ts  T^{\tss t}_1(u).
\een
Hence \eqref{transp} is an anti-homomorphism.
Finally, for \eqref{inverse} observe that the relation
\ben
R(u-v)\ts T^{-1}_2(v)\ts  T^{-1}_1(u)  =
T^{-1}_1(u)\ts  T^{-1}_2(v)\ts R(u-v)
\een
is equivalent to \eqref{ternary}.
Note now that the mappings
\eqref{sign} and \eqref{transp} are involutive and so these two
anti-homomorphisms are bijective.

The bijectivity of the anti-homomorphism $\Sr$ of $\Y(\gl_N)$
defined by \eqref{inverse} follows from the bijectivity
of its square $\Sr^{2}$ which is computed at the end of this
section.
\epf

The anti-automorphisms \eqref{sign} and \eqref{transp} are involutive
and commute with each other.
Their composition is an involutive automorphism
of $\Y(\gl_N)$ such that
\beql{sin}
T(u)\mapsto T^{\tss t}(-u)\,.
\eeq
This automorphism of the algebra $\Y(\gl_N)$
will play an important role in Section  \ref{sec:provec}.
However, the
anti-automorphism \eqref{inverse} is {\em not\/} involutive unless $N=1\,$.
This is the antipodal map $\Sr$
of the Hopf algebra $\Y(\gl_N)$, see Section~\ref{sec:hopf} below. 

To compute the square of the anti-homomorphism \eqref{inverse}
consider $N\times N$ matrix obtained from $T^{-1}(u)$ by transposition. 
Let us denote this new matrix by $T^{\tss\sharp}(u)\,$.
Accordingly, the $ij$-th entry
of this matrix will be denoted by $T^{\tss\sharp}_{ij}(u)\,$.
This entry is a formal power series in $u^{-1}$
with coefficients from the algebra $\Y(\gl_N)\,$. By definition,
\begin{equation}
\label{SrT}
\Sr:T_{ij}(u)\mapsto T^{\tss\sharp}_{ji}(u)\,.
\end{equation}
Our computation of the image of $T_{ij}(u)$ relative to $\Sr^{2}$
is based on the next lemma. 

\ble
\label{zu}
There is a formal power series $Z(u)$ in $u^{-1}$ with the coefficients from
the centre of the algebra\/ $\Y(\gl_N)$ and with 
the leading term $1$ such that for all $i$~and~$j$
\beql{zetav}
\sum_{k=1}^N\,T_{ki}(u+N)\,T^{\tss\sharp}_{kj}(u)=\delta_{ij}\,Z(u)\,.
\eeq
\ele

\bpf
Let us multiply both sides of the relation \eqref{ternary}
by $T_2^{-1}(v)$ on the left and right and then apply transposition relative
to the second copy of $\End\CC^N$. We get
$$
R^{\tss t}(u-v)\,T^{\tss\sharp}_2(v)\ts T_1(u)=
T_1(u)\,T^{\tss\sharp}_2(v)\ts R^{\tss t}(u-v)\,.
$$
Multiplying both sides of this result on the left and right 
$R^{\tss t}(u-v)^{-1}$ we get
\beql{trater}
R^{\tss t}(u-v)^{-1}\,
T_1(u)\,T^{\tss\sharp}_2(v)=
T^{\tss\sharp}_2(v)\ts T_1(u)\,
R^{\tss t}(u-v)^{-1}\,.
\eeq
Multiplying the latter equality by $u-v-N$ and then setting
$u=v+N$ we get
\beql{qresi}
Q\,T_1(v+N)\,T^{\tss\sharp}_2(v)=T_2^{\tss\sharp}(v)\,T_1(v+N)\,Q\,,
\eeq
see \eqref{rtinv}.
Because the operator $Q$ is one-dimensional,
either side of \eqref{qresi}
must be equal to $Q$ times a certain power series in $v^{-1}$
with the coefficients from $\Y(\gl_N)\,$.  Denote
this series by $Z(v)\,$.
By applying the left hand side of \eqref{qresi} to the basis vector
$e_i\ot e_j$ we obtain the required equality \eqref{zetav}.

It is immediate from \eqref{tiju} and \eqref{zetav} that the leading
term of series $Z(v)$ is $1\,$. Let us prove that all the coefficients
of this series are central in $\Y(\gl_N)\,$.
We will work with the algebra \eqref{multitp}
where $m=3\,$. By using the relations \eqref{ternary} and \eqref{trater},
\begin{gather*}
R^{\tss t}_{13}(u-v)^{-1}\,R_{12}(u-v-N)\,
T_1(u)\,T_2(v+N)\,T^{\tss\sharp}_3(v)=
\\
R^{\tss t}_{13}(u-v)^{-1}\,
T_2(v+N)\,T_1(u)\,T^{\tss\sharp}_3(v)\,
R_{12}(u-v-N)=
\\
T_2(v+N)\,T^{\tss\sharp}_3(v)\,T_1(u)\,
R^{\tss t}_{13}(u-v)^{-1}\,R_{12}(u-v-N)\,.
\end{gather*}
Note that by using the expressions \eqref{rmatrix} and \eqref{rtinv}
we obtain the equality
$$
Q_{23}\,R^{\tss t}_{13}(u-v)^{-1}\,R_{12}(u-v-N)=
Q_{23}\,(\ts1-(u-v-N)^{-2}\ts).
$$
So multiplying the first and third lines of
previous display by $Q_{23}$ on the left gives
$$
(\ts1-(u-v-N)^{-2}\ts)\,
T_1(u)\,Z(v)\,Q_{23}\,=\,
Q_{23}\,Z(v)\,T_1(u)\,
(\ts1-(u-v-N)^{-2}\ts)
$$
where we also used \eqref{qresi}.
The last display shows that any generator $T_{ij}^{(r)}$
commutes with every coefficient of the series $Z(v)\,$.
\epf

It follows from \eqref{tiju} and \eqref{zetav} 
that the coefficient of the series
$Z(u)$ at $u^{-1}$ is zero. In Section \ref{sec:cen} we will show that 
the coefficients of $Z(u)$ at $u^{-2},u^{-3},\ldots$ are free generators
of the centre of the algebra $\Y(\gl_N)\,$. Hence we will again use
Lemma~\ref{zu}.

\bpr
\label{ss}
The square of the map\/ $\Sr$ is the automorphism of\/
$\Y(\gl_N)$ given~by
\ben
\Sr^{\tss2}:T(u)\mapsto Z(u)^{-1}\ts T(u+N)\,.
\een
\epr

\bpf
Let us apply the anti-homomorphism $\Sr$ to both sides of
the identity
\ben
\sum_{k=1}^N\,T_{jk}(u)\,T^{\tss\sharp}_{ik}(u)=\delta_{ij}\,.
\een
Using \eqref{SrT} we get
$$
\sum_{k=1}^N\,\Sr^2(T_{ki}(u))\ts T^{\tss\sharp}_{kj}(u)
=\delta_{ij}\,.
$$
Comparing this with \eqref{zetav} we conclude that
$
\,\Sr^2(T_{ki}(u))=Z(u)^{-1}\,T_{ki}(u+N)\,.
$
\epf


\section{Hopf algebra structure}
\label{sec:hopf}

A {\it coalgebra\/} over the field $\CC$ is a complex vector space $\Ar$
equipped with a linear map $\Delta:\Ar\to\Ar\ot\Ar$ called the
{\it comultiplication\/}, and another linear map
$\ve:\Ar\to\CC$ called the {\it counit\/},
such that the following three diagrams are commutative:

\ben
\begin{CD}
\Ar 
@>{\De}>>
\Ar\ot\Ar
\\
@V{\De}VV
@VV{\De\ot\text{id}}V
\\
\Ar\ot\Ar
@>>{\text{id}\ot\De}>               
\Ar\ot\Ar\ot\Ar
\end{CD}
\vspace{8pt}
\een

\noindent
which gives the {\it coassociativity\/} axiom of 
the comultiplication $\Delta\,$, and

\ben
\begin{CD}
\Ar
@>{\De}>>
\Ar\ot\Ar
\\
@V{\text{id}}VV
@VV{\ve\ot\text{id}}V
\\
\Ar
@>>{\cong}>
\CC\ot\Ar 
\end{CD}
\qquad\qquad
\begin{CD}
\Ar
@>{\De}>>
\Ar\ot\Ar
\\
@V{\text{id}}VV
@VV{\text{id}\ot\ve}V
\\
\Ar
@>>{\cong}>
\Ar\ot\CC 
\end{CD}
\vspace{8pt}
\een

A {\it bialgebra\/} over $\CC$ is a complex associative unital
algebra $\Ar$ equipped with a coalgebra structure,
such that $\De$ and $\ve$ are algebra homomorphisms.
In particular, then $\De(1)=1\ot1$ and $\ve(1)=1$.
A bialgebra $\Ar$ is called a {\it Hopf algebra},
if it is also equiped with an anti-automorphism
$\Sr:\Ar\to\Ar$ called the {\it antipode\/},
such that another two diagrams are commutative:

\ben
\begin{CD}
\Ar
@>{\de\,\ve}>>
\Ar
\\
@V{\De}VV
@AA{\mu}A
\\
\Ar\ot\Ar
@>>{\Sr\ot\text{id}}>
\Ar\ot\Ar
\end{CD}
\qquad\qquad
\begin{CD}
\Ar
@>{\de\,\ve}>>
\Ar
\\
@V{\De}VV
@AA{\mu}A
\\
\Ar\ot\Ar
@>>{\text{id}\ot\Sr}>
\Ar\ot\Ar
\end{CD}
\vspace{8pt}
\een

\noindent
Here $\mu:\Ar\ot\Ar\to\Ar$ is the algebra multiplication
and $\de:\CC\to\Ar$ is the unit map of the algebra $\Ar$,
that is $\de(c)=c\cdot1$ for any $c\in\CC$.

\bpr\label{thm:hopf}
The Yangian $\Y(\gl_N)$ is a Hopf algebra with comultiplication
\beql{Delta}
\Delta: T_{ij}(u)\mapsto
\sum_{k=1}^N
\,
T_{ik}(u)\ot T_{kj}(u),
\eeq
the antipode \eqref{inverse}
and the counit $\ve:T(u)\mapsto 1$.
\epr

\bpf
We start by verifying the axiom that
$\Delta:\Y(\gl_N)\to\Y(\gl_N)\ot\Y(\gl_N)$
is an algebra homomorphism. We shall slightly generalize the notation used 
in Section~\ref{sec:matrix}. Let $m$ and $n$ be positive integers.
Introduce the algebra
\beql{mnprod}
(\End\CC^N)^{\otimes m}\otimes
\Y(\gl_N)^{\ot n}.
\eeq
For all $a\in\{1,\dots,m\}$ and $b\in\{1,\dots,n\}$ consider
the formal power series in $u^{-1}$ with the coefficients in this algebra,
\ben
T_{a[b]}(u)=\sum_{i,j=1}^N\,
1^{\ot(a-1)}\ot e_{ij}\ot 1^{\otimes(m-a)}
\ot
1^{\ot(b-1)}\ot T_{ij}(u)\ot 1^{\otimes(n-b)}.
\een
The definition of $\Delta$ can now be written in a matrix form,
\beql{deltamatrix}
\Delta: T(u)\mapsto T_{[1]}(u)\ts T_{[2]}(u)
\eeq
where $T_{[b]}(u)$ is an abbreviation for
the series $T_{1[b]}(u)$ with the coefficients from
the algebra \eqref{mnprod} where $m=1$ and $n=2$.
We need to show
that $\Delta (T(u))$ obeys \eqref{ternary}:
\begin{gather*}
R(u-v)\,T_{1[1]}(u)\,T_{1[2]}(u)\,T_{2[1]}(v)\,T_{2[2]}(v)
\,=\,
\\
T_{2[1]}(v)\,T_{2[2]}(v)\,T_{1[1]}(u)\,T_{1[2]}(u)\,R(u-v).
\end{gather*}
Here $m=n=2$, and $R(u-v)$ is identified with $R(u-v)\ot1\ot1$.
But this relation
is implied by the relation \eqref{ternary},
and by the observation that the elements $T_{1[2]}(u)$ and
$T_{2[1]}(v)$ commute, as well as
the elements $T_{1[1]}(u)$ and $T_{2[2]}(v)$ do.

Our $\Sr$ is an anti-automorphism 
relative to multiplication due to Proposition~\ref{prop:anti}. 
Since $\De$ is a homomorphism of algebras, the definition
\eqref{deltamatrix} implies that 
$$
\Delta:T^{-1}(u)\mapsto T_{[2]}^{-1}(u)\ts T_{[1]}^{-1}(u)\,.
$$
Therefore $\Sr$ is also an anti-automorphism relative to comultiplication.
The other~two axioms involving $\Sr$ are
readily verified since
\ben
(\Sr\ot\text{id})\,\Delta: T(u)\mapsto T^{-1}_{[1]}(u)\ts
T^{}_{[2]}(u)
\een
and
\ben
(\text{id}\ot\Sr)\,\Delta: T(u)\mapsto T^{}_{[1]}(u)\ts T^{-1}_{[2]}(u)
\een
so that subsequent application of $\mu$ yields the identity
matrix in both the cases. 
\epf

We have $\ve\,\bigl(T_{ij}^{(r)}\bigr)=0$ for $r\geqslant1\,$.
By expanding the formal power series in $u^{-1}$ in \eqref{Delta}
we obtain a more explicit definition of the comultiplication $\Delta$
on $\Y(\gl_N)\,$,
\beql{Deltamatelem}
\Delta\bigl(T_{ij}^{(r)}\bigr)=
T_{ij}^{(r)}\ot1+1\ot T_{ij}^{(r)}+
\,
\sum_{k=1}^N
\sum_{s=1}^{r-1}
\,
T_{ik}^{(s)}\ot T_{kj}^{(r-s)}\,.
\eeq
Hence this comultiplication is not cocommutative unless $N=1\,$.

\bpr
\label{lz}
For the series $Z(u)$ defined above we have
$$
\De:Z(u)\mapsto Z(u)\ot Z(u)\,.
$$
\epr

\bpf
The square $\Sr^{\tss2}$ of
of the antipodal map is a coalgebra automorphism. 
Hence the images of $T(u)$ relative to the compositions
$\De\,\Sr^{\tss2}$ and $(\Sr^{\tss2}\ot\Sr^{\tss2})\,\De$
are the same. By Proposition \ref{ss} these images are respectively equal to
$$
\De\,(\tss Z(u)^{-1}\,T(u+N))=
\De\,(\tss Z(u)^{-1})\,(\tss T(u+N)\ot T(u+N))
$$
and
$$
(\Sr^{\tss2}\ot\Sr^{\tss2})(\tss T(u)\ot T(u))=
(\tss Z(u)^{-1}\,T(u+N))\ot(\tss Z(u)^{-1}\,T(u+N))\,.
$$
Here we identify $Z(u)^{-1}$ with the series
$1\ot Z(u)^{-1}$ which takes its coefficients from 
$\End\CC^N\otimes\Y(\gl_N)$ and use the homomorphism
property of $\De\,$.
By dividing the right hand sides of above two equalities by 
$T(u+N)\ot T(u+N)$ and equating the results 
\ben
\De:Z(u)^{-1}\mapsto Z(u)^{-1}\ot Z(u)^{-1}\,.
\qedhere
\een
\epf


\section{Two filtrations on the Yangian}
\label{sec:filtr}

There are two natural ascending filtrations on
the associative algebra $\Y(\gl_N)\,$. The first one is defined by
\ben
\deg T_{ij}^{(r)}=r\,.
\een
For any $r\geqslant 1$ we will denote by $\Tb_{ij}^{(r)}$ the 
image of the generator $T_{ij}^{(r)}$ 
in the degree $r$ component of the corresponding graded algebra
$\gr\Y(\gl_N)\,$.
It is immediate from the defining relations \eqref{defequiv}
that all these images pairwise commute.
In Section \ref{sec:pbw} we will prove that these images 
are also algebraically independent.


Now introduce another filtration on $\Y(\gl_N)$ by setting for $r\geqslant1$
\beql{secondfiltr}
\degpr T_{ij}^{(r)}=r-1\,.
\eeq
Let $\grpr\Y(\gl_N)$ be the corresponding graded algebra.
Let $\Tt_{ij}^{\ts(r)}$ be the image of $T_{ij}^{(r)}$
in the component of $\grpr\Y(\gl_N)$ of the degree $r-1\,$.

The graded algebra $\grpr\Y(\gl_N)$
inherits from $\Y(\gl_N)$ the Hopf algebra structure. Namely,
by using \eqref{Deltamatelem} for any $r\geqslant1$ we get
\begin{gather}
\label{gradedhopf1}
\De\bigl(\,\Tt_{ij}^{\ts(r)}\bigr)=
\Tt_{ij}^{\ts(r)}\ot1+1\ot\Tt_{ij}^{\ts(r)},
\\
\label{gradedhopf2}
\ve\bigl(\,\Tt_{ij}^{\ts(r)}\bigr)=0
\fand
{\rm S}\bigl(\,\Tt_{ij}^{\ts(r)}\bigr)=
-\tss\Tt_{ij}^{\ts(r)}.
\end{gather}

For any Lie algebra $\g$ over the field $\CC$
consider the universal enveloping algebra $\U(\g)$.
There is a natural Hopf algebra structure on $\U(\g)$.
The comultiplication $\Delta$,
the counit $\ve$ and  the antipode ${\rm S}$ on $\U(\g)$
are defined by setting for $X\in\g$
\begin{gather}
\label{standardhopf1}
\Delta(X)=X\ot1+1\ot X\,,
\\[2pt]
\label{standardhopf2}
\ve(X)=0
\fand
{\rm S}(X)=-X\,.
\end{gather}
In the next proposition $\g$ is the polynomial current Lie algebra
$\gl_N[z]\cong\gl_N\ot\CC[z]\,$. 
The latter Lie algebra is naturally graded
by degrees of the indeterminate $z\,$. 

\bpr
\label{prop:iso2}
The graded Hopf algebra\/ $\grpr\Y(\gl_N)$ is isomorphic to
$\U(\gl_N[z])$.
\epr

\bpf
Using the defining relations \eqref{defequiv} we get
\ben
[\,\Tt_{ij}^{\ts(r)}, \Tt_{kl}^{\ts(s)}] =
\de_{kj}\,\Tt_{il}^{\ts(r+s-1)}-
\de_{il}\,\Tt_{kj}^{\ts(r+s-1)}.
\een
Hence the assignments 
\beql{asshom}
E_{ij}\tss z^{r-1}\mapsto\Tt_{ij}^{\ts(r)}
\quad\text{for}\quad
r\geqslant1
\eeq
define a surjective homomorphism
\beql{surhom}
\U(\gl_N[z])\to\grpr\Y(\gl_N)
\eeq
of graded associative algebras.
At the end of Section \ref{sec:pbw} we will show that
the kernel of this homomorphism is trivial.
Hence comparing the definitions \eqref{gradedhopf1},\eqref{gradedhopf2}
with the general definitions 
\eqref{standardhopf1},\eqref{standardhopf2}
completes the proof of Proposition \ref{prop:iso2}.
\epf



\section{Vector and covector representations}
\label{sec:provec}

We shall often use the matrix $T(u)$ to describe
homomorphisms from $\Y(\gl_N)$ to other algebras.
Namely, let $\Ar$ be any unital associative algebra over the
field $\CC\tss$. Let $X(u)$ be the $N\times N$ matrix whose
$ij$-entry is any formal power series $X_{ij}(u)$ in $u^{-1}$
with the leading term $\de_{ij}$ and all coefficients
from the algebra $\Ar\,$. If $\al:\Y(\gl_N)\to\Ar$ is
any homomorphism, then
the assignment
\beql{assign}
\al:T(u)\mapsto X(u)
\eeq
means that every coefficient of the series $T_{ij}(u)$
gets mapped to the corresponding coefficient of the series $X_{ij}(u)$
for all indices $i,j=1,\dots,N$.
If we regard $T(u)$ as a series in $u$ with the coefficients
from the algebra $\End\CC^N\ot\Y(\gl_N)$ then,
more formally, we may write
\ben
{\id}\ot\al:T(u)\mapsto X(u)
\een
instead of \eqref{assign}. Here
\ben
X(u)=\sum_{i,j=1}^N e_{ij}\ot X_{ij}(u),
\een
is regarded as a series in $u$ with
coefficients from the algebra $\End\CC^N\ot\Ar\tss$;
cf.\ \eqref{tutens}.

Setting $\Ar=\End\CC^N$ and $X(u)=R(u)$ above,
we can define a homomorphism $\Y(\gl_N)\to\End\CC^N$
by the assignment $T(u)\mapsto R(u)$.
To prove the homomorphism property
by using the matrix form \eqref{ternary} of the defining relations
of the algebra $\Y(\gl_N)\tss$,
we have to check the equality of rational functions in $u$ and $v$
with values in the algebra $(\End \CC^N)^{\ot3}$,
\ben
R_{12}(u-v)\ts R_{13}(u)\ts R_{23}(v)=
R_{23}(v)\ts R_{13}(u)\ts R_{12}(u-v)\tss.
\een
But this equality is just another form of \eqref{ybe}.
In other words, the assignment $T(u)\mapsto R(u)$ defines a representation
of $\Y(\gl_N)$ on the vector space $\CC^N$. Here
\ben
T_{ij}(u)\mapsto\de_{ij}-e_{ji}\,u^{-1}
\een
by \eqref{pmatrix} and \eqref{rmatrix}.
Note that
this representation of the algebra $\Y(\gl_N)$ can also be obtained by
pulling the defining representation $E_{ij}\mapsto e_{ij}$
of the Lie algebra $\gl_N$ back through
the homomorphism \eqref{pin}.
This remark justifies the definition \eqref{pin}.

By pulling the defining representation $E_{ij}\mapsto e_{ij}$
of the Lie algebra
$\gl_N$ back through the homomorphism \eqref{eval}, we get
the representation of $\Y(\gl_N)$ such that
\ben
T_{ij}(u)\mapsto\de_{ij}+e_{ij}\,u^{-1}\tss.
\een
Hence this representation can be described by the assignment
$T(u)\mapsto R^{\tss t}(-u)\,$. Observe that the
representations $T(u)\mapsto R(u)$ and $T(u)\mapsto R^{\tss t}(-u)$
differ by the involutive automorphism \eqref{sin}
of the algebra $\Y(\gl_N)\tss$.

By pulling the representation $T(u)\mapsto R(u)$
back through the automorphism \eqref{shift} of $\Y(\gl_N)$
for any $c\in\CC\tss$,
we get the representation of $\Y(\gl_N)$ on
the vector space $\CC^N$, such that $T(u)\mapsto R(u-c)$.
It is called a \textit{vector representation\/} of
$\Y(\gl_N)\tss$, and is denoted by $\rho_c\tss$. Thus
\ben
\rho_c:
T_{ij}(u)\mapsto\de_{ij}-e_{ji}\,(u-c)^{-1}
\een
or equvalently,
\beql{3.66}
\rho_c:
T_{ij}^{(r)}\mapsto-\tss c^{\tss r-1}\tss e_{ji}
\,\quad\text{for any}\quad
r\geqslant1\tss.
\eeq

By pulling the representation $T(u)\mapsto R^{\tss t}(-u)$
back through the automorphism \eqref{shift},
we get the representation of $\Y(\gl_N)$ on
$\CC^N$, such that $T(u)\mapsto R^{\tss t}(c-u)$.
It is called a \textit{covector representation\/} of
$\Y(\gl_N)\tss$, and is denoted by $\si_c\tss$. Thus
\ben
\si_c:
T_{ij}(u)\mapsto\de_{ij}+e_{ij}\,(u-c)^{-1}
\een
or equvalently,
\beql{3.666666}
\si_c:
T_{ij}^{(r)}\mapsto c^{\tss r-1}\tss e_{ij}
\,\quad\text{for any}\quad
r\geqslant1\tss.
\eeq

In Section \ref{sec:filtr} we introduced
an ascending filtration on algebra $\Y(\gl_N)$
such that any generator $T_{ij}^{(r)}$ of $\Y(\gl_N)$ has the
degree $r-1\,$. We denoted the corresponding graded algebra by
$\grpr\Y(\gl_N)$ and defined a surjective homomorphism \eqref{surhom}
by \eqref{asshom}.

Under this homomorphism the element
$T_{ij}^{(r)}$ of $\Y(\gl_N)\tss$,
or rather its image $\Tt^{\tss(r)}_{ij}$ in
$\grpr\Y(\gl_N)\,$, corresponds to the generator
$E_{ij}\tss z^{\,r-1}$ of $\U(\gl_N[z])$.
One can define a representation $\st_c$
of the algebra $\U(\gl_N[z])$ on the vector space
$\CC^N$ by
\beql{bareval}
\st_c:
E_{ij}\tss z^{r-1}\mapsto c^{\tss r-1}\tss e_{ij}
\,\quad\text{for any}\quad
r\geqslant1\tss,
\eeq
so that
\ben
\st_c(\tss E_{ij}\tss z^{r-1}\tss)=\si_c(\tss T_{ij}^{(r)}\tss)\,.
\een
The representation $\st_c$ is an example of an
\text{evaluation representation} of $\U(\gl_N[z])$,
see the general definition in Section \ref{sec:eval} below.


\section{Evaluation representations}
\label{sec:eval}

For any Lie algebra
$\agot$ over $\CC$ consider the corresponding polynomial current Lie
algebra $\agot\tss[z]=\agot\ot\CC[z]\,$.
Let $\om$ be any representation of $\agot$
on the vector space $\CC^N$, and $c$
be any complex number.
Then one can define a representation of $\agot\tss[z]$ by
\ben
X\tss z^s\mapsto c^{\tss s}\om\tss(X)
\,\quad\text{for any}\quad
s\geqslant0\tss.
\een
This is the \textit{evaluation representation\/} of the Lie
algebra $\agot\tss[z]\tss$, corresponding to $\om$ at the point
$z=c$ of the complex plane $\CC\tss$.
When $\agot=\gl_N$ and $\om$ is the defining representation
of the Lie algebra $\gl_N$ on $\CC^N$,
we obtain $\st_c$ in this way.

We will need the following general property of evaluation
representations.
For any $c_1,\dots,c_n\in\CC$ let us denote by $\om_{c_1\dots c_n}$ the
tensor product of the evaluation representations of
the Lie algebra $\agot\tss[z]$
corresponding to $\om$ at the points $c_1,\dots,c_n\tss$.
We extend the representation $\om_{c_1\dots c_n}$
to the universal enveloping algebra $\U(\agot\tss[z])\tss$.

\ble
\label{lem:L2.1}
Suppose that the Lie algebra $\agot$ is finite-dimensional,
and $\om$ is its faithful representation.
Let the parameters $c_1,\dots,c_n$ and integer $n\geqslant0$
vary. Then the intersection in $\U(\agot\tss[z])$ of the
kernels of all representations $\om_{c_1\dots\tss c_n}$ is trivial.
\ele

\bpf
Using the faithful representation $\om$ of the Lie algebra $\agot\tss$,
we can identify $\agot\tss[z]$ with a subalgebra of the Lie algebra
$\gl_N[z]$.
Hence it suffices to consider the case when $\agot$ is the
Lie algebra $\gl_N\tss$, and $\om:\gl_N\to\End\CC^N$ is the
defining representation. Let us assume that this is the case.
Then $\om_c=\st_c$ as we have already observed.

Let us now choose any basis $X_1,\dots,X_{N^2}$ of $\gl_N$
such that its first vector $X_1$ is 
\ben
I=E_{11}+\dots+E_{NN}\tss.
\een
To distinguish between the algebras $\U(\gl_N)$ and $\End\CC^N$,
the operators on $\CC^N$
corresponding to the elements $X_1,\dots,X_{N^2}\in\gl_N$
will be denoted by $x_1,\dots,x_{N^2}$ respectively.
Note that then $x_1$ is the identity operator $1$.

The elements $X_{a}\tss z^{\ts s}$ 
with $a=1,\dots,N^2$ and $s=0\ts,1\ts,2\ts,\ts\dots$
constitute a basis of $\gl_N\tss[z]\,$.
Choose any total ordering of this basis which ends with the 
infinite~sequence
$$
\dots\,,\ts X_1\tss z^{\ts 2},\tss X_1\ts z\,,\ts X_1\,.
$$
Take any finite linear combination $L$ of the products
\beql{2.3.1}
(X_{a_1}z^{\ts s_1})\dots(X_{a_m}z^{\ts s_m})\in\U(\gl_N[z])
\eeq
with 
$$
L_{\,a_1\dots\ts a_m}^{\,s_1\dots\ts s_m}\in\CC
$$
being the coprrespondinf coefficients. 
The number $m$ of factors in
\eqref{2.3.1} may~vary.
Assume that the factors in each product \eqref{2.3.1} are
arranged according to our chosen ordering of the basis of $\gl_N\tss[z]\,$.
Due to the commutation relations in $\U(\gl_N\tss[z])$ 
we may assume it without any loss of generality.
Suppose that $\om_{c_1\dots\ts c_n}(L)=0$ for all $n$
and $c_1,\dots,c_n\in\CC\,$. We need to prove that $L=0\,$.

For each product \eqref{2.3.1} there is a number $p$ such that 
the indices
$a_1,\dots,a_{\ts p}>1$ but $a_{\ts p+1},\dots,a_m=1\,$.
This is due to our ordering of the basis of $\gl_N\tss[z]\,$.
The numbers $p$ for different products
\eqref{2.3.1} may differ, and we do not exclude the case $p=0$ here.
Let $h$ be the maximum of the numbers $p$ in our linear combination $L\,$.

Suppose that $n\geqslant h\,$.
Let $\omega_{\ts h}$ be the symmetrisation map of
the tensor product $(\gl_N\tss[z])^{\ot h}$ normalised so that 
$\omega_{\ts h}^{\ts2}=h\ts!\,\omega_{\ts h}\,$.
Let $V$ be the subspace
of $(\End\CC^N)^{\ot n}$ spanned by the vectors
$x_{\ts b_{\ts1}}\ot\dots\ot x_{\ts b_{\ts n}}$
where at least one of the 
indices $b_{\ts1},\dots,b_{\,h}$ is $1\,$.
If $h=0$ then this subspace is assumed to be zero.
Applying the homomorphism 
$\om_{c_1\dots\ts c_n}$ to a product \eqref{2.3.1} 
with $p=h$ gives
\begin{gather}
\label{2.3.2}
(\ts\om_{c_1}\ot\dots\ot\om_{c_{h}}\ts)\,
(\,\omega_{\ts h}\ts
(\ts X_{\ts a_1}z^{\ts s_1}
\ot\dots\ot 
X_{\ts a_{h}}z^{\ts s_h}\ts)\ts)
\ot1^{\ts\ot\ts(n-h)}\ \times
\\
\nonumber
\prod_{k=h+1}^m
(\ts c_1^{\ts s_k}+\dots+c_n^{\ts s_k}\ts)
\end{gather}
modulo the subspace $V$. Applying
$\om_{c_1\dots\ts c_n}$ to a product \eqref{2.3.1} 
with $p<h$ gives an element of $V$.
But a linear combination of the 
expressions \eqref{2.3.2} belongs to $V$
only if this combination is zero.

For each product \eqref{2.3.1} with $p=h$
there is a certain number $l\geqslant h$ such that 
$s_{\ts h+1}\geqslant\dots\geqslant s_l>0$ 
but $s_{\ts l+1},\dots,s_m=0\,$.
This is due to our ordering of the basis of $\gl_N\tss[z]\,$.
Then \eqref{2.3.2} equals $n^{\ts m-l}$ times
\begin{gather}
\label{2.3.3}
(\ts\om_{c_1}\ot\dots\ot\om_{c_{h}}\ts)\,
(\,\omega_{\ts h}\ts
(\ts X_{\ts a_1}z^{\ts s_1}
\ot\dots\ot 
X_{\ts a_{h}}z^{\ts s_h}\ts)\ts)
\ot1^{\ts\ot\ts(n-h)}\ \times
\\
\nonumber
\prod_{k=h+1}^l
(\ts c_1^{\ts s_k}+\dots+c_n^{\ts s_k}\ts)\,.
\end{gather}
Let $g$ be the maximum of the numbers $l$
for all products \eqref{2.3.1} with $p=h\,$.
 
Suppose that $n\geqslant g\,$. Consider the pairs of sequences
$a_1,\ldots,a_m$ and $s_1,\dots,s_m$
showing in $L$ in any product \eqref{2.3.1} with $p=h\,$.
For every such pair there is  some $l\in\{h,\ldots g\}\,$.
Then $a_{\ts h+1},\dots,a_{m}=1$ and $s_{\ts l+1},\dots,s_m=0\,$.
Take all different pairs of sequences 
$a_1,\ldots,a_h$ and $s_1,\dots,s_l$ arising in this way.
The expressions \eqref{2.3.3} corresponding to the latter
pairs of sequences
are linearly independent
as polynomials in $c_1,\dots,c_n$ 
with values in $(\End\CC^N)^{\ot n}\,$.
This is again is due to our ordering of the basis of $\gl_N\tss[z]\,$.
Here we also use the observation that in \eqref{2.3.3} the image of
$\om_{c_1}\ot\dots\ot\om_{c_{h}}$
does not depend on the parameters $c_{h+1},\dots,c_n$
whereas the product over $k=h+1,\dots,l$ in \eqref{2.3.3}
depends on these parameters when $l>h\,$.

Therefore if $\om_{c_1\dots\ts c_n}(L)=0$ for a certain $n\geqslant g$
and for all $c_1,\dots,c_n\in\CC$ then 
\begin{equation*}
\sum_{m=l}^\infty\ 
L_{\,a_1\dots a_h\,1\dots\ts 1}^{\,s_1\dots\,s_l\,0\dots\ts0}\ 
n^{\ts m-l}=0\,.
\end{equation*}
In the last displayed sum there are exactly $m$ lower indices
and also $m$ upper indices in the coefficient
$L_{\,a_1\dots a_h\,1\dots\ts 1}^{\,s_1\dots\,s_l\,0\dots\ts0}\,$. 
By letting the number $n$ vary we now prove that all these 
coefficients vanish. 
\epf


\section{Poincar\'e--Birkhoff--Witt theorem}
\label{sec:pbw}

Let us now make use of the bialgebra structure on $\Y(\gl_N)\tss$.
For any $c_1,\dots,c_n\in\CC$ take the tensor product
of the vector representations $\rho_{c_1},\dots,\rho_{c_n}$
of $\Y(\gl_N)\tss$. We get a representation
\ben
\rho_{c_1\dots\tss c_n}:\Y(\gl_N)\to(\End\CC^N)^{\ot n}.
\een
If $n=0$, the representation $\rho_{c_1\dots\tss c_n}$
is understood as the counit homomorphism $\ve:\Y(\gl_N)\to\CC\tss$.
Using the matrix form
\eqref{deltamatrix} of the definition of the
comultiplication on $\Y(\gl_N)\tss$, we see that
\ben
\id\ot\rho_{c_1\dots\tss c_n}:
T(u)\mapsto
R_{12}(u-c_1)\dots R_{1,n+1}(u-c_n)\tss.
\een
Here we apply the convention made in
the beginning of Section \ref{sec:provec}
to the algebra $\Ar=(\End\CC^N)^{\ot n}$
and to the homomorphism $\al=\rho_{c_1\dots\tss c_n}$.

The tensor product of the covector representations
$\si_{c_1},\dots,\si_{c_n}$ will be denoted by
$\si_{c_1\dots\tss c_n}$. By using the matrix form
\eqref{deltamatrix} of the definition of the
comultiplication on $\Y(\gl_N)$ again, we see that
\ben
\id\ot\si_{c_1\dots\tss c_n}:
T(u)\mapsto
R^{\tss t}_{12}(c_1-u)\dots R^{\tss t}_{1,n+1}(c_n-u)\tss.
\een
By using Lemma~\ref{lem:L2.1}, we will now prove the following proposition.

\bpr
\label{thm:P2.2}
Let the parameters $c_1,\dots,c_n\in\CC$ and the integer $n\geqslant0$
vary. Then the intersection of the
kernels of all representations $\si_{c_1\dots\tss c_n}$ is trivial.
\epr

\bpf
Take any finite linear combination $A$ of the products
\ben
T_{i_1j_1}^{(r_1)}\dots T_{i_mj_m}^{(r_m)}\in\Y(\gl_N)
\een
with certain coefficients
\ben
A_{\tss i_1j_1\dots i_mj_m}^{\,r_1\dots r_m}\in\CC
\een
where the indices $r_1\tss,\dots,r_m\geqslant1$ and the number
$m\geqslant0$ may vary, as well as the indices
$i_1,j_1,\dots,i_m,j_m\tss$. Suppose that $A\neq0$ as
an element of $\Y(\gl_N)\tss$.

The algebra $\Y(\gl_N)$ comes with an ascending filtration
such that $T_{ij}^{(r)}$ has the degree $r-1\,$. Let $d$ be the
degree of $A$ rtelative to this filtration. Let $B$ be the image of $A$
in the degree $d$ component of the graded algebra $\grpr\Y(\gl_N)\tss$.
Then $B\neq0$.

We can also assume that
\ben
A_{\tss i_1j_1\dots i_mj_m}^{\,r_1\dots r_m}=0
\ \quad\text{if}\quad
r_1+\dots+r_m>d+m\tss.
\een
Let $C$ be the
sum of the elements of the algebra $\U(\gl_N[z])$,
\ben
\sum_{r_1+\dots+r_m=d+m}
A_{\tss i_1j_1\dots i_mj_m}^{\,r_1\ldots r_m}\,
(E_{\tss i_1j_1}\tss z^{\tss r_1-1})
\dots
(E_{\tss i_mj_m}\tss z^{\tss r_m-1})\tss.
\een
The image of $C$ under the homomorphism \eqref{surhom} equals $B$.
In particular, $C\neq0$.

Consider the image of $A$ under the representation
$\si_{c_1\dots\tss c_n}$. This image
depends on $c_1,\dots,c_n$ polynomially.
The degree of this polynomial does not exceed $d$
by the definition \eqref{3.666666}.
Let $D$ be the sum of the terms of degree $d$
of this polynomial.

Now equip the tensor product $\Y(\gl_N)^{\ot n}$
with the ascending filtration
where the degree is the sum of the degrees on the tensor factors.
Then under the $n$-fold comultiplication
$\Y(\gl_N)\to\Y(\gl_N)^{\ot n}$
\ben
T_{ij}^{(r)}\mapsto
\sum_{b=1}^n\,1^{\ot(b-1)}\ot T_{ij}^{(r)}\ot1^{\ot(n-b)}
\ \ \text{plus terms of degree less than}\ \ r-1\,,
\een
see \eqref{Deltamatelem}. But under the
$n$-fold comultiplication
$\U(\gl_N[z])\to\U(\gl_N[z])^{\ot n}$,
\ben
E_{ij}\tss z^{r-1}\mapsto
\sum_{b=1}^n\,1^{\ot(b-1)}\ot(E_{ij}\tss z^{r-1})\ot1^{\ot(n-b)}\tss.
\een
The definitions \eqref{3.666666} and \eqref{bareval} now imply
that the sum $D\in(\End\CC^N)^{\ot n}$
coincides with the image of the sum $C\in\U(\gl_N[z])$
under the tensor product of the evaluation representations
$\st_{c_1},\dots,\st_{c_n}$.
Since $C\neq0$, using Lemma ~\ref{lem:L2.1}
we can choose $n$ and $c_1,\dots,c_n$ so that
$D\neq0\tss$. Then $\si_{c_1\dots c_n}(A)\neq0\tss$
by the definition of $D$.
\epf

\bpr
\label{thm:P2.3}
Let the parameters $c_1,\dots,c_n\in\CC$ and the integer $n\geqslant0$
vary. Then the intersection of the
kernels of all representations $\rho_{c_1\dots\tss c_n}$ is trivial.
\epr

The proof of Proposition \ref{thm:P2.3} is similar
to that of Proposition \ref{thm:P2.2} and is omitted.
We will now prove the injectivity of homomorphism \eqref{surhom}
by modifying the logic of our proof of Proposition \ref{thm:P2.2}.
Take any finite linear combination $C$ of the products
\ben
(E_{\tss i_1j_1}\tss z^{\tss r_1-1})
\dots
(E_{\tss i_mj_m}\tss z^{\tss r_m-1})
\in\U(\gl_N[z])
\een
with certain coefficients
\ben
C_{\tss i_1j_1\dots i_mj_m}^{\,r_1\dots r_m}\in\CC
\een
where the indices $r_1\tss,\dots,r_m\geqslant1$ and the number
$m\geqslant0$ may vary, as well as the indices
$i_1,j_1,\dots,i_m,j_m\tss$. Suppose that $C\neq0$ as
an element of $\U(\gl_N[z])\tss$.

The algebra $\U(\gl_N[z])$ is graded so that
for any integer $s\geqslant0\tss$,
the generator $E_{ij}\tss z^s$ has the degree $s$.
The homomorphism \eqref{surhom} preserves the degree.
Without loss of generality suppose that the element $C$ is
homogeneous of degree $d\,$, that is
\ben
C_{\tss i_1j_1\dots i_mj_m}^{\,r_1\dots r_m}=0
\ \quad\text{if}\quad
r_1+\dots+r_m\neq d+m\tss.
\een
Now define the element $A\in\Y(\gl_N)$ as the sum
\ben
\sum_{r_1+\dots+r_m=d+m}
C_{\tss i_1j_1\dots i_mj_m}^{\,r_1\ldots r_m}\,
T_{i_1j_1}^{(r_1)}\dots T_{i_mj_m}^{(r_m)}\tss.
\een
Let $B$ be the image of $A$
in the $d\tss$-th component of the graded algebra $\grpr\Y(\gl_N)\tss$.
The element $B$ coincides with the image of $C$
under the homomorhism \eqref{surhom}.

Now let $D\in(\End\CC^N)^{\ot n}$ be the image of $C$
under the tensor product of the evaluation representations
$\st_{c_1},\dots,\st_{c_n}$.
The image of $A$ under the representation
$\si_{c_1\dots\tss c_n}$
depends on $c_1,\dots,c_n$ polynomially.
The degree of this polynomial does not exceed $d$
by \eqref{3.666666}.
The sum of the terms of degree $d$
of this polynomial equals $D$, see the proof of Proposition~\ref{thm:P2.2}.
Since $C\neq0$, using Lemma ~\ref{lem:L2.1}
we can choose $n$ and $c_1,\dots c_n$ so that
$D\neq0\tss$. Then $\degpr A=d$.
Indeed, if $\degpr A<d$ then the degree of
the polynomial $\si_{c_1,...,c_n}(A)$ would be also less
then $d$. This would contradict to the non-vanishing of 
$D$. By the definition of the element $B\in\grpr\Y(\gl_N)\tss$,
the equality $\degpr A=d$ means that $B\neq0$.
So the homomorphism \eqref{surhom} is injective.

Let us now invoke the classical
Poincar\'e--Birkhoff--Witt theorem for the universal
enveloping algebras of Lie algebras \cite[Section~2.1]{d:ae}.
By applying this theorem to the Lie algebra $\gl_N[z]$ we now obtain its
analogue for the Yangian $\Y(\gl_N)\,$.

\bth
\label{thm:pbw}
Given an arbitrary linear ordering of the set of
generators $T^{(r)}_{ij}$ with $r\geqslant 1\,$,
any element of the algebra $\Y(\gl_N)$ can be uniquely
written as a linear combination of ordered monomials in these generators.
\eth


\bco\label{cor:polalg}
The graded algebra 
$\gr \Y(\gl_N)$ is the algebra of polynomials in the generators
$\Tb_{ij}^{\ts(r)}$ with $r\geqslant 1\,$.
\eco


\section{Centre of the Yangian}
\label{sec:cen}

Let $\agot$ be any Lie algebra over the field $\CC\,$.
Consider the corresponding polynomial current Lie algebra
$\agot\tss[z]\,$.
In the proof of Theorem \ref{thm:zcentre} we will use a
general property of the universal enveloping algebra $\U(\agot\tss[z])\,$.
It is stated as the lemma below.
 
\ble
\label{lem:centcur}
Suppose that the Lie algebra $\agot$ is finite-dimensional and has
the trivial centre. Then the centre of the algebra $\U(\agot\tss[z])$
is also trivial, that is equal to\/ $\CC\,$.
\ele

\bpf
Consider adjoint action of the Lie algebra $\agot\tss[z]$ 
on its symmetric algebra. It suffices to prove that the space
of invariants of this action is trivial.

Let $A$ be any element of the symmetric algebra of $\agot\tss[z]$
invariant under the adjoint action. Let $M=\dim\agot$.
Choose any basis $X_1,\ldots,X_M$ of $\agot$ and let
\ben
[\,X_p\,,X_q\,]=\sum_{r=1}^M c_{pq}^{\tss r}\,X_r
\een
where $c_{pq}^{\tss r}\in\CC\,$. 
Let $L$ be the minimal non-negative integer such that
\ben
A\,=\sum_{d_1,\ldots,d_M}A_{\,d_1\ldots\tss d_M}\,
(X_1\tss  z^{\tss L})^{d_1}\ldots(X_M\tss z^{\tss L})^{d_M}
\een
where 
$d_1,\ldots,d_M$ range over 
non-negative integers and
$A_{\,d_1\ldots d_M}$ is a polynomial in the basis elements
$X_p\tss z^{\tss s}$ of $\agot\tss[z]$ with $1\leqslant p\leqslant M$ and
$0\leqslant s<L$ only.
We have 
$$
\text{ad}\tss(X_p\tss z)(A)=0 
$$
for every index $p=1,\dots,M$.
The component of the left hand side of this equation
that involves the basis elements of $\agot\tss[z]$
of the form $X_r\tss z^{\tss L+1}$
must be zero. Thus
\ben
\sum_{d_1,\ldots,d_M}
A_{\,d_1\ldots\tss d_M}
\sum_{q,r=1}^M
c_{pq}^{\tss r}\,d_q\,
(X_1\tss z^{\tss L})^{d_1}
\ldots
(X_q\tss z^{\tss L})^{d_q-1}
\ldots
(X_M\tss z^{\tss L})^{d_M}
\,X_r\tss z^{\tss L+1}=0\,.
\een
Taking here the coefficient of $X_r\tss z^{\tss L+1}$ we obtain that
\ben
\sum_{d_1,\ldots,d_M}
A_{\,d_1\ldots\tss d_M}
\,\sum_{q=1}^M\,
c_{pq}^{\tss r}\,d_q\,
(X_1\tss z^{\tss L})^{d_1}
\ldots
(X_q\tss z^{\tss L})^{d_q-1}
\ldots
(X_M\tss z^{\tss L})^{d_M}=0\,.
\een
If follows that for any non-negative integers 
$d_1^{\,\prime},\ldots,d_M^{\,\prime}$ we have
\beql{adprim}
\sum_{q=1}^M\,
A_{\,d_1^{\tss\prime}\ldots\tss d^{\tss\prime}_q+1\ldots\tss d^{\tss\prime}_M}
c_{pq}^{\tss r}\,(d^{\tss\prime}_q+1)=0
\quad\text{where}\quad
p,r=1,\ldots,M\,.
\eeq

Let us now fix $d_1^{\tss\prime},\ldots,d_M^{\tss\prime}$ 
and observe that the elements
$X_q^{\tss\prime}=(d_q^{\tss\prime}+1)\ts X_q$
with $q=1,\ldots,M$
also form a basis of $\agot\,$. Since the centre of $\agot$ is trivial,
the system 
\ben
[\,X_p\,,\sum_{q=1}^M\, a_q\tss X_q^{\tss\prime}\,]=0
\quad\text{where}\quad
p=1,\ldots,M
\een
of linear equations on 
$a_1,\ldots,a_M\in\CC$
has only trivial solution. 
It can be written as
\ben
\sum_{q=1}^M\,a_q\,c_{pq}^{\tss r}\,(d_q^{\tss\prime}+1)=0
\quad\text{where}\quad
p,r=1,\ldots,M\,.
\een
Hence by comparing \eqref{adprim} with the latter system
we obtain that 
$A_{\,d_1^{\tss\prime}\ldots\tss 
d_q^{\tss\prime}+1\ldots\tss d_M^{\tss\prime}}=0$ for every
$q=1,\ldots,M\,$. It now follows that $A\in\CC\,$,
and $L=0$ in particular.
\epf

Now consider the series $Z(u)$ defined by \eqref{zetav}.
For any $r\geqslant 1$ let $Z^{(r)}$ be the coefficient
of this series at $u^{-r}$. Just before stating Proposition \ref{ss}
we noted that $Z^{(1)}=0\,$. Hence
$$
Z(u)=1+Z^{(2)}u^{-2}+Z^{(3)}u^{-3}+\ldots\,.
$$

\bpr
\label{zim}
For any $r\geqslant2$ the element 
$Z^{(r)}\in\Y(\gl_N)$ has the degree $r-2$
relative to the filtration \eqref{secondfiltr}. Its image
in the graded algebra $\grpr\Y(\gl_N)$ is equal to
$$
(\tss1-r\tss)\,\sum_{i=1}^N\,\Tt^{\tss(r-1)}_{ii}\,.
$$
\epr

\bpf
Let us expand the factor $T_{ki}(u+N)$ appearing in the definition 
\eqref{zetav} as a formal power series in $u^{-1}\,$. The result 
has the form
$$
T_{ki}(u)+N\,\dot{T}_{ki}(u)+X_{ki}(u)
$$ 
where
$$
\dot{T}_{ki}(u)=
-\,T^{(1)}_{ki}u^{-2}
-2\,T^{(2)}_{ki}u^{-3}-\ldots
$$
is the formal derivative of the series $T_{ki}(u)$ and
$$
X_{ki}(u)=X^{(3)}_{ki}u^{-3}+X^{(4)}_{ki}u^{-4}+\ldots
$$
is a series with coefficients $X^{(r)}_{ki}\in\Y(\gl_N)$
such that $\degpr X_{ki}^{(r)}=r-3$ for $r\geqslant3\,$.
By setting $i=j$ in \eqref{zetav} and summing
over $i=1,\ldots,N$ we now get the equality
\begin{equation}
\label{nzu}
N+
\sum_{i,k=1}^N\,
(\tss N\,\dot{T}_{ki}(u)+X_{ki}(u))\,T^{\tss\sharp}_{ki}(u)=N\tss Z(u)\,.
\end{equation}
Here we used the definition
of the matrix $T^{\tss\sharp}(u)$ as the transposed inverse of $T(u)\,$.

The leading term of the series $T^{\tss\sharp}_{ki}(u)$
is $\de_{ik}$ while for any $r\geqslant1$ the coefficient of this series
at $u^{-r}$ has the degree $r-1$ relative to \eqref{secondfiltr}.
It follows that modulo lower degree elements,
for any $r\geqslant1$ the coefficient at $u^{-r}$
of the series at the left hand side of \eqref{nzu} equals
$$
N\,\sum_{i,k=1}^N\,(1-r)\,T^{\tss(r-1)}_{ki}\,\de_{ik}
=
N\,\sum_{i=1}^N\,(1-r)\,T^{\tss(r-1)}_{ii}\,.
$$
Hence Proposition \ref{zim} follows from 
\eqref{nzu}. Also we see once again that $Z^{(1)}=0\,$.
\epf

\bth
\label{thm:zcentre}
The coefficients $Z^{(2)},Z^{(3)},\ldots$ of the series $Z(u)$
are free generators of the centre of the associative algebra $\Y(\gl_N)\,$.
\eth

\bpf
Let us apply Lemma \ref{lem:centcur} to the special linear Lie algebra 
$\agot=\sll_N\,$.
Since the centre of the universal enveloping algebra 
$\U(\tss\sll_N\tss[z]\tss)$ is trivial, the 
decomposition 
$$
\gl_N\tss[z]=\sll_N\tss[z]\,\oplus\,\CC\tss[z]\,\sum_{i=1}^N\,E_{ii}
$$
of Lie algebras
implies that the centre of $\U(\tss\gl_N\tss[z]\tss)$
is generated by the elements
\beql{cengen}
\sum_{i=1}^N\,E_{ii}\tss z^{r-1}
\eeq
where $r\geqslant1\,$. Moreover these generators
are free due to
the Poincar\'e--Birkhoff--Witt theorem \cite[Section~2.1]{d:ae}
applied to the Lie algebra $\gl_N[z]\,$.

Under the isomorphism \eqref{surhom},
the elements \eqref{cengen} go respectively to the elements
$$
\sum_{i=1}^N\,\Tt^{\tss(r)}_{ii}
$$
where again $r\geqslant1\,$. Therefore the latter elements are free
generators of the centre of the algebra $\grpr\Y(\gl_N)\,$, see
Proposition \ref{prop:iso2}.
On the other hand, we have already proved that 
the elements $Z^{(2)},Z^{(3)},\ldots$ of the algebra $\Y(\gl_N)$
belong to its centre, see Lemma \ref{zu}. Hence Theorem 
\ref{thm:zcentre} follows from Proposition \ref{zim}.
\epf


\section{Dual Yangian}
\label{sec:dualyang}

The {\it dual Yangian\/} for
the Lie algebra $\gl_N\,$,
denoted by $\Y^\ast(\gl_N)$, is an
associative unital algebra
over the field $\CC$
with a countable set of generators
$$
T_{ij}^{(-1)},\ T_{ij}^{(-2)},\dots
\quad\text{where}\quad
i,j=1,\dots,N\,.
$$
To write down the defining relations for these generators,
introduce the 
series
\beql{3.10}
T^{\tss\ast}_{ij}(v)=
\de_{ij}+T_{ij}^{\ts(-1)}+
T_{ij}^{\ts(-2)}\ts v+T_{ij}^{\ts(-3)}v^2+\ldots
\in\Y^{\ast}(\gl_N)[[v]]\,.
\eeq
The reason for separating the term $\de_{ij}$ in \eqref{3.10}
will become apparent in the next section. Now
combine all the series \eqref{3.10} into the single element
\beql{tutensast}
T^{\tss\ast}(v)=
\sum_{i,j=1}^N T^{\tss\ast}_{ij}(v)\ot e_{ij}\in\Y^{\ast}(\gl_N)[[v]]\ot
\End\CC^N\,.
\eeq

We will write the defining relations of the algebra $\Y^\ast(\gl_N)$
first in their matrix form, to be compared with \eqref{ternary}.
For any positive integer $n$, consider the algebra
\beql{multitpast}
\Y^\ast(\gl_N)\ot(\End\CC^N)^{\ot n}.
\eeq
For any index $b\in\{1,\dots,n\}$
introduce the formal power series in the variable $v$
with the coefficients from the algebra \eqref{multitpast},
\beql{tbast}
T^{\tss\ast}_b(v)=\sum_{i,j=1}^N T^{\tss\ast}_{ij}(v)\ot
1^{\ot(b-1)}\ot e_{ij}\ot 1^{(n-b)}.
\eeq
Here the 
belongs to the $b$-th copy of $\End\CC^N$.
Setting $n=2$ and identifying $R(u-v)$ with $1\ot R(u-v)\,$,
the defining relations of  $\Y^\ast(\gl_N)$ can be written as
\beql{3.13}
T^{\tss\ast}_1(u)\ts T^{\tss\ast}_2(v)\ts R(u-v)=
R(u-v)\ts T^{\tss\ast}_2(v)\ts T^{\tss\ast}_1(u)\,.
\eeq

The relation \eqref{3.13} is
equivalent to the collection of relations
\ben
(u-v)\ts [\tss T^{\tss\ast}_{ij}(u),T^{\tss\ast}_{kl}(v)\tss]=
T^{\tss\ast}_{il}(u)\,T^{\tss\ast}_{kj}(v)-
T^{\tss\ast}_{il}(v)\,T^{\tss\ast}_{kj}(u)
\een
for all $i,j,k,l=1,\dots,N$.
We omit the proof of the equivalence, because
it is similar to the proof of Proposition~\ref{prop:ternary}.
The last displayed relation can be rewritten as
\ben
[\tss T^{\tss\ast}_{ij}(u),T^{\tss\ast}_{kl}(v)\tss]=
\sum_{p=0}^{\infty}\,
u^{-p-1}v^p\,
\Big(
T^{\tss\ast}_{il}(u)\,T^{\tss\ast}_{kj}(v)-
T^{\tss\ast}_{il}(v)\,T^{\tss\ast}_{kj}(u)
\Big).
\een 
Expanding here the series in $u,v$ and
equating the coefficients at $u^{r-1}v^{s-1}$ we get
\begin{gather}
\label{3.12}
[\,T^{\tss(-r)}_{ij}, T^{\tss(-s)}_{kl}]
\,=\,
\de_{kj}\,T^{\tss(-r-s)}_{il}-\de_{il}\,T^{\tss(-r-s)}_{kj}\,+
\\
\nonumber
\sum_{b=1}^{s}
\Big(
T^{\tss(b-r-s-1)}_{il}\,T^{\tss(-b)}_{kj}
-
T^{\tss(-b)}_{il}\,T^{\tss(b-r-s-1)}_{kj}
\Big).
\end{gather}
The proof of next proposition is similar to that of
Proposition~\ref{thm:hopf} and is omitted.

\bpr
\label{prop:dualhopf}
The dual Yangian $\Y^\ast(\gl_N)$ is a bialgebra
over the field\/ $\CC$ with the counit defined
$\ve:T^{\tss\ast}(v)\mapsto1$ and the comultiplication defined by
\beql{3.17}
\Delta: T^{\tss\ast}_{ij}(v)\mapsto\sum_{k=1}^N
T^{\tss\ast}_{ik}(v)\ot T^{\tss\ast}_{kj}(v)\,.
\eeq
\epr

Expanding the power series in $v$ in \eqref{3.17}
and using the axiom $\Delta(1)=1\ot1$,
we get a more explicit definition of the comultiplication
on the dual Yangian $\Y^\ast(\gl_N)$,
\beql{Deltamadual}
\Delta\bigl(T_{ij}^{(-r)}\bigr)=
T_{ij}^{(-r)}\ot1+1\ot T_{ij}^{(-r)}+
\sum_{k=1}^N
\sum_{s=1}^r
\,
T_{ik}^{(-s)}\ot T_{kj}^{(s-r-1)}
\eeq
for $r\geqslant1$; cf. \eqref{Deltamatelem}.
Since $\ve(1)=1$,
for every $r\geqslant1$ we get $\ve\bigl(T_{ij}^{(-r)}\bigr)=0\tss$.

The dual Yangian $\Y^\ast(\gl_N)$ is a bialgebra but
{\it not\/} a Hopf algebra.
The antipodal map $\Sr$ is defined only for a completion
$\Y^\circ(\gl_N)$ of $\Y^\ast(\gl_N)$ such that the element
\ben
T^{\tss\ast}(0)\in\Y^\circ(\gl_N)\ot\End\CC^N
\een is invertible.
Then 
\ben
T^{\tss\ast}(v)\in\Y^\circ(\gl_N)[[v]]\ot\End\CC^N
\een
is also invertible, and the antipode $\Sr$ is defined by
mapping $T^{\tss\ast}(v)$ to its inverse. This inverse will be 
denoted by $T^{\tss\natural}(v)$. 
It will be used again in the end of Section~\ref{sec:dpbw}.

In order to construct such a completion, let us
equip the algebra $\Y^\ast(\gl_N)$ with a
descending filtration, defined by assigning to the generator
$T_{ij}^{(-r)}$ the degree $r$ for any $r\geqslant1$.
Then $\Y^\circ(\gl_N)$ is defined
as the formal completion of 
$\Y^\ast(\gl_N)$ relative to this descending
filtration.
The algebra $\Y^\circ(\gl_N)\ot\End\CC^N$ 
contains the inverse~of 
\ben
T^{\tss\ast}(0)=1\ot1+\sum_{i,j=1}^N T_{ij}^{(-1)}\ot e_{ij}\,.
\een

We extend the comultiplication $\Delta$ on $\Y^\ast(\gl_N)$ to 
$\Y^\circ(\gl_N)\,$, and also denote this extension by $\Delta\,$. The image
$\Delta(\Y^\circ(\gl_N))$ lies in the formal completion of
the algebra $\Y^\ast(\gl_N)\ot\Y^\ast(\gl_N)$
with respect to the descending filtration,
defined by assigning to the element $T_{ij}^{(-r)}\ot T_{kl}^{(-s)}$
the degree $r+s$. Indeed, the image $\Delta\bigl(T_{ij}^{(-r)}\bigr)$
in $\Y^\ast(\gl_N)\ot\Y^\ast(\gl_N)$
is a sum of elements of degrees $r$ and $r+1$
by \eqref{Deltamadual}.

The kernel of the counit homomorphism
$\ve:\Y^\ast(\gl_N)\to\CC$ consists of all the elements
which of positive degree relative to the filtration,
see Proposition~\ref{prop:dualhopf}.
Therefore $\ve$ extends to the algebra $\Y^\circ(\gl_N)$.
This extension is the counit map for the Hopf algebra
$\Y^\circ(\gl_N)$, it will be also denoted by $\ve$.

For any $c\in\CC$
the assignment $T^{\tss\ast}(v)\mapsto T^{\tss\ast}(v+c)$ determines
an automorphism of the algebra $\Y^\circ(\gl_N)$.
This follows from the relations \eqref{3.13},
cf.\ Proposition~\ref{prop:auto}. But for $c\neq0$
this automorphism does not preserve the subset
$\Y^\ast(\gl_N)\subset\Y^\circ(\gl_N)$, and therefore does not
determine an automorphism of $\Y^\ast(\gl_N)\,$.

To find the square of the antipodal map $\Sr$ of the Hopf algebra
$\Y^\circ(\gl_N)$ let $T^{\tss\flat}(v)$ be the result of applying
to the inverse 
of \eqref{tutensast} the 
transposition in 
$\End\CC^N\,$.~Write
$$
T^{\tss\flat}(v)=
\sum_{i,j=1}^N T^{\tss\flat}_{ij}(v)\ot e_{ij}\in\Y^{\circ}(\gl_N)[[v]]\ot
\End\CC^N
$$
so that 
$$
\Sr:T_{ij}^{\tss\ast}(v)\mapsto T^{\tss\flat}_{ji}(v)\,.
\vspace{4pt}
$$
The proof of the next lemma is similar 
to that of Lemma \ref{zu} and is omitted here.

\ble
\label{zc}
There is a formal power series $Z^{\tss\circ}(v)$ in $v$ 
with coefficients from the centre of the algebra $\Y^\circ(\gl_N)$
such that for all indices $i$~and~$j$
$$
\sum_{k=1}^N\,T_{ki}^{\tss\ast}(v-N)\,T^{\tss\flat}_{kj}(v)=
\delta_{ij}\,Z^{\tss\circ}(v)\,.
$$
\ele

In general, the coefficients of the series $Z^{\tss\circ}(v)$
do not belong to the dual Yangian $\Y^\ast(\gl_N)\,$. 
However, the proposition below
can be derived from Lemma \ref{zc} just as 
Proposition \ref{ss} was derived from Lemma \ref{zu}.
Hence we again omit the proof.

\bpr
\label{sc}
The square of the map\/ $\Sr$ is the automorphism of\/
$\Y^\circ(\gl_N)$ 
\ben
\Sr^{\tss2}:T^{\tss\ast}(v)\mapsto 
Z^{\tss\circ}(v)^{-1}\ts T^{\tss\ast}(v-N)\,.
\een
\epr

The latter result follows just
as Proposition \ref{lz} followed from Proposition~\ref{ss}.

\bpr
\label{zzcc}
For the series $Z^{\tss\circ}(v)$ defined above we have
$$
\De:Z^{\tss\circ}(v)\mapsto Z^{\tss\circ}(v)\ot Z^{\tss\circ}(v)\,.
$$
\epr

The completion
$\Y^\circ(\gl_N)$ of the filtered algebra $\Y^\ast(\gl_N)$
can be described more explicitly. At the end of Section \ref{sec:nondeg}
we will show
that the vector space $\Y^\ast(\gl_N)$ has a basis parameterized
by all multisets of triples $(r_1,i_1,j_1),\dots,(r_m,i_m,j_m)$ where
\ben
r_1,\dots,r_m\in\{1,2,\dots\}
\fand
i_1,j_1,\dots,i_m,j_m\in\{1,\dots,N\}
\een
while $m=0,1,2,\dots\ts\tss$.
The corresponding basis vector in $\Y^\ast(\gl_N)$
is the monomial
\beql{dualbasis}
T_{i_1j_1}^{\ts(-r_1)}
\dots\,
T_{i_mj_m}^{\ts(-r_m)}.
\eeq
The ordering of the factors in this monomial can be chosen
arbitrarily.
Choose any linear ordering of the basis monomials.
For any positive integer $r$, there is only a
finite number of the basis monomials \eqref{dualbasis}
such that $r_1+\dots+r_m\leqslant r$.
This means that when the index of the basis monomial
\eqref{dualbasis} in any chosen linear ordering increases, then
the filtration degree \eqref{dualbasis}
\ben
r_1+\dots+r_m\to\infty\,.
\een
Therefore the vector space $\Y^\circ(\gl_N)$ consists of
all infinite linear combinations of the basis monomials
\eqref{dualbasis}, with the  coefficients from the field $\CC\tss$.


\section{Canonical pairing}
\label{sec:pairing}

There is a canonical bilinear pairing
\beql{pairing}
\langle\,\ts,\,\rangle:\Y(\gl_N)\times\Y^\ast(\gl_N)\to\CC\tss.
\eeq
We shall describe the corresponding linear map
$\beta:\Y(\gl_N)\otimes\Y^\ast(\gl_N)\to\CC$.
It will be defined so that
for all integers $m,n\geqslant0$ the linear map
\ben
(\End\CC^N)^{\ot m}\ot\Y(\gl_N)\ot\Y^\ast(\gl_N)\ot(\End\CC^N)^{\ot n}
\to
(\End\CC^N)^{\ot(m+n)}
\een
given by $\text{id}\ot\beta\ot\text{id}$, will send
\beql{3.7777777}
T_1(u_1)\dots T_m(u_m)\ot T^{\tss\ast}_{1}(v_1)\dots T^{\tss\ast}_n(v_n)
\,\mapsto\!\!\!
\prod_{1\leqslant a\leqslant m}^\rightarrow\,
\prod_{1\leqslant b\leqslant n}^\rightarrow R_{a,b+m}(u_a-v_b).
\eeq
Here $u_1,\dots,u_m,v_1,\dots,v_n$ are independent
variables. The coefficients of the series
\beql{series}
T_1(u_1),\dots,T_m(u_m)
\fand
T^{\tss\ast}_{1}(v_1),\dots,T^{\tss\ast}_n(v_n)
\eeq
belong to the algebras \eqref{multitp} and \eqref{multitpast},
respectively.

Note that the series in $u_1,\dots,u_m$ and $v_1,\dots,v_n$
at the left hand side of 
\eqref{3.7777777}
satisfy certain relations, implied by the defining relations
of the algebras $\Y(\gl_N)$ and $\Y^\ast(\gl_N)$.
The following proposition guarantees that
the pairing
is well-defined.

\bpr\label{prop:agreement}
The assignment \eqref{3.7777777} agrees
with 
relations \eqref{ternary} and \eqref{3.13}.
\epr

\bpf
This follows from
the Yang-Baxter equation \eqref{ybe}.
For instance, let us consider the case when $m=2$ and $n=1$. Here
we have to check that the series
\ben
\big(\ts R(u_1-u_2)\,T_1(u_1)\,T_2(u_2)\tss\big)\ot T^{\tss\ast}(v)
\een
and
\ben
\big(\ts T_2(u_2)\,T_1(u_1)\,R(u_1-u_2)\tss\big)\ot T^{\tss\ast}(v)
\een
with the coefficients in the algebra
\ben
(\End\CC^N)^{\ot2}\ot\Y(\gl_N)\ot\Y^\ast(\gl_N)\ot\End\CC^N,
\een
have the same images in
under the
map $\text{id}\ot\beta\ot\text{id}$.
These images are series with the coefficients in
$(\End\CC^N)^{\ot3}$.
Note that the second element
can be rewritten as
\ben
\big(\ts P\,T_1(u_2)\,T_2(u_1)\,P\,R(u_1-u_2)\tss\big)\ot T^{\tss\ast}(v)
\een
By the definition \eqref{3.7777777}, the images of the
two elements are respectively
\ben
R_{12}(u_1-u_2)\ts R_{13}(u_1-v)\ts R_{23}(u_2-v)
\een
and
\begin{gather*}
P_{12}\ts R_{13}(u_2-v)\ts R_{23}(u_1-v)\ts P_{12}\ts R_{12}(u_1-u_2)
\\
\,=\,
R_{23}(u_2-v)\ts R_{13}(u_1-v)\ts R_{12}(u_1-u_2).
\end{gather*}
The equality of two images is now evident due to \eqref{ybe}.
Using \eqref{ybe} repeatedly, one can prove
Proposition \ref{prop:agreement}
for any $m,n\geqslant0$.
\epf

Let us show that the assignments \eqref{3.7777777}
for all $m,n=0,1,2,\tss\dots$
determine the values of the bilinear pairing \eqref{pairing} uniquely.
When $m=n=0$, we get from \eqref{3.7777777}
the equality $\langle\,1\,,1\,\rangle=1\,$.
By choosing $m=1$ and $n=0$, we obtain from
\eqref{3.7777777} that
$\langle\ts T_{ij}^{(r)},1\ts\rangle=0$ for any $r\geqslant1\,$.
When $m=0$ and $n=1$, we obtain that
$\langle\ts 1,T_{ij}^{(-s)}\ts\rangle=0$ for any $s\geqslant1\,$.
In both cases, we had to use the equality
$\langle\,1\,,1\,\rangle=1$ obtained above.

Now suppose that $m,n\geqslant1$. To determine the pairing values
\beql{pvalue}
\bigl\langle\,
T_{i_1j_1}^{(r_1)}\dots T_{i_mj_m}^{(r_m)}\,,\,
T_{k_1l_1}^{(-s_1)}\dots T_{k_nl_n}^{(-s_n)}
\bigl\rangle
\eeq
for any indices
\ben
r_1,\dots,r_m,s_1,\dots,s_n\in\{1,2,\dots\}
\een
and
\ben
i_1,j_1,\dots,i_m,j_m,k_1,l_1,\dots,k_n,l_n\in\{1,\dots,N\},
\een
the product of the rational functions
$R_{a,b+m}(u_a-v_b)$ on the right hand side of \eqref{3.7777777}
should be expanded as power series in the variables
$u_1^{-1},\dots,u_m^{-1},v_{1},\dots,v_n\,$.
The series \eqref{series}
should be then also expanded.

Note that although the coefficient of $v^0$ in the series
\eqref{3.10} is a sum of two terms,
$\de_{ij}$ and $T_{ij}^{\ts(-1)}$,
the pairing value \eqref{pvalue} can be still determined
by \eqref{3.7777777} for any indices $s_1,\dots,s_n\geqslant1$
by using the induction on $n$. Namely, if some of the indices
$s_1,\dots,s_n$ are equal to $1$, the value \eqref{pvalue}
can be determined by \eqref{3.7777777}, using
the values \eqref{pvalue} with
$n$ replaced by $0,\dots,n-1$.

Consider the case $m=n=1$ in more detail. Then the map
$\text{id}\ot\beta\ot\text{id}$
maps
\ben
T(u)\ot T^{\tss\ast}(v)=
\sum_{i,j,k,l=1}^N
e_{ij}\ot
\biggl(\,
\sum_{r=0}^\infty\,
T_{ij}^{(r)}\,u^{-r}
\biggr)
\ot
\biggl(
\de_{kl}+\sum_{s=1}^\infty\,T_{kl}^{(-s)}v^{s-1}
\biggr)
\ot e_{kl}
\een
to the series
\beql{rexp}
R(u-v)=
1\ot1-\sum_{i,j=1}^N\,\sum_{r=1}^\infty\,
u^{-r}v^{r-1}\,e_{ij}\ot e_{ji}\,;
\eeq
see \eqref{rmatrix} and \eqref{3.10}.
Using the equality
$\langle\ts T_{ij}^{(r)},1\ts\rangle=\de_{0r}$ for $r\geqslant0$, we get
\ben
\bigl\langle\,
T_{ij}^{(r)},\tss
T_{kl}^{(-s)}
\bigr\rangle
=-\tss\de_{rs}\,\de_{il}\,\de_{jk}
\quad\text{for}\quad
r,s\geqslant1.
\een

More explicitly the value \eqref{pvalue} will be determined 
in the course of the proof of the next lemma. This lemma
describes a basic property of the
bilinear pairing \eqref{pairing}.
It is valid for any integers $m,n\geqslant0$.

\ble\label{lem:L4.1}
If\/ $r_1+\dots+r_m<s_1+\dots+s_n$
then the value \eqref{pvalue} is zero.
\ele

\bpf
First suppose that $s_1,\dots,s_n\geqslant2$.
Then by the definition
of the pairing \eqref{3.7777777},
the value \eqref{pvalue} is the coefficient of
\beql{3.91}
e_{i_1j_1}\ot\dots\ot e_{i_mj_m}
\ot
e_{k_1l_1}\ot\dots\ot e_{k_nl_n}
\cdot\,
u_1^{-r_1}\dots u_m^{-r_m}
\,
v_1^{s_1-1}\dots v_n^{s_n-1}
\eeq
in the expansion of the product in
$(\End\CC^N)^{\ot(m+n)}
[[u_1^{-1},\dots,u_m^{-1},v_1,\dots,v_n]]$
\ben
\prod_{1\leqslant a\leqslant m}^\rightarrow
\,
\prod_{1\leqslant b\leqslant n}^\rightarrow
R_{a,b+m}(u_a-v_b)
=
\prod_{1\leqslant a\leqslant m}^\rightarrow
\,
\prod_{1\leqslant b\leqslant n}^\rightarrow\!
\bigl(\,
1-\sum_{r=1}^\infty\,u_a^{-r}v_b^{r-1}P_{a,b+m}
\,\bigr)\,.
\een
If the coefficient of \eqref{3.91}
is non-zero in this expansion then clearly we have the inequality
$r_1+\dots+r_m\geqslant s_1+\dots+s_n\,$.

Now suppose that some of the numbers $s_1,\dots,s_n$ are equal to 1.
Without loss of generality we will assume that $s_1,\dots,s_d\geqslant2$
and $s_{d+1},\dots,s_n=1$ for some $d<n$.
Rewrite the product
at the right hand side of the definition \eqref{3.7777777} as
\ben
\prod_{1\leqslant b\leqslant d}^\rightarrow
\,
\prod_{1\leqslant a\leqslant m}^\rightarrow
R_{a,b+m}(u_a-v_b)
\ \cdot\!
\prod_{d<b\leqslant n}^\rightarrow
\,
\prod_{1\leqslant a\leqslant m}^\rightarrow
R_{a,b+m}(u_a-v_b)
\,.
\een
By definition, the coefficient of \eqref{3.91} in the
expansion of this product equals
\ben
\bigl\langle\,
T_{i_1j_1}^{(r_1)}\dots T_{i_mj_m}^{(r_m)}\,,\,
T_{k_1l_1}^{(-s_1)}\dots T_{k_dl_d}^{(-s_d)}
\bigl(\de_{k_{d+1}l_{d+1}}\!+\,T_{k_{d+1}l_{d+1}}^{(-1)}\bigr)
\dots
\bigl(\de_{k_nl_n}\!+\,T_{k_nl_n}^{(-1)}\bigr)
\bigl\rangle
\,.
\een
The value \eqref{pvalue} is then
the coefficient of \eqref{3.91} in the expansion of the product
\beql{3.925}
\prod_{1\leqslant b\leqslant d}^\rightarrow
\,
\prod_{1\leqslant a\leqslant m}^\rightarrow
\bigl(\,
1-\sum_{r=1}^\infty\,u_a^{-r}v_b^{r-1}P_{a,b+m}
\,\bigr)\ \times
\eeq
\ben
\prod_{d<b\leqslant n}^\rightarrow
\Bigl(\,
\prod_{1\leqslant a\leqslant m}^\rightarrow
\bigl(\,
1-\sum_{r=1}^\infty\,u_a^{-r}v_b^{r-1}P_{a,b+m}
\,\bigr)
-1\Bigr)\,.
\een
If that coefficient here is non-zero, then
$r_1+\dots+r_m\geqslant s_1+\dots+s_d+n-d\,$.
\epf


\section{Non-degeneracy of the pairing}
\label{sec:nondeg}

In Section \ref{sec:dualyang}
we equipped the algebra $\Y^\ast(\gl_N)$
with a descending filtration. Now consider
the corresponding graded algebra $\gr\Y^\ast(\gl_N)$.
Its component of degree $s$ will be
denoted by $\text{gr}_s\Y^\ast(\gl_N)$. For any $s\geqslant1$
denote by $\Tt_{ij}^{\ts(-s)}$
the image of $T_{ij}^{(-s)}$ in 
$\text{gr}_s\Y^\ast(\gl_N)$.
By \eqref{3.12} we immediately~get
\ble
\label{lem:grel}
In the graded algebra\/ $\gr\Y^\ast(\gl_N)$,
for any $r,s\geqslant1$ we have
\ben
[\,\Tt^{\ts(-r)}_{ij},\Tt^{\ts(-s)}_{kl}]
=\tss
\de_{kj}\,\Tt^{\ts(-r-s)}_{il}
-\tss
\de_{il}\,\Tt^{\ts(-r-s)}_{kj}
.
\een
\ele

In Section \ref{sec:filtr} we equipped the algebra $\Y(\gl_N)$ with
an ascending filtration, such that the corresponding graded algebra
$\gr\Y(\gl_N)$ is commutative.
Its subspace of all elements of degree $s$ will be
denoted by $\text{gr}_s\Y(\gl_N)$.
Keeping to the notation of Section \ref{sec:filtr},
for any $s\geqslant1$ let $\Tb_{ij}^{\ts(s)}$
be the image of the generator $T_{ij}^{(s)}$ in
$\text{gr}_s\Y(\gl_N)$.

We can define a bilinear pairing
\beql{3.93}
\langle\,\ts,\ts\rangle:\,\,\gr\Y(\gl_N)\times\gr\Y^\ast(\gl_N)\to\CC
\eeq
by making its value
\beql{gpvalue}
\bigl\langle\,
\Tb_{i_1j_1}^{\ts(r_1)}
\dots
\Tb_{i_mj_m}^{\ts(r_m)}
\,,\,
\Tt_{k_1l_1}^{\ts(-s_1)}
\dots
\Tt_{k_nl_n}^{\ts(-s_n)}
\,\bigr\rangle
\eeq
equal to \eqref{pvalue} if
$
r_1+\ldots+r_m=s_1+\ldots+s_n
$
and by making it equal to zero otherwise.
Here 
$r_1,\dots,r_m,s_1,\dots,s_n\geqslant1$
and 
$m,n\geqslant0\,$.
The indices $i_1,j_1,\dots,i_m,j_m$ and $k_1,l_1,\dots,k_n,l_n$ may be
arbitrary. This definition is self-consistent. Namely, if
\beql{self-cons}
r_1+\dots+r_m=s_1+\ldots+s_n=s
\eeq
for some $s\geqslant1$, then by Lemma~\ref{lem:L4.1} we have
$$
\bigl\langle
T_{i_1j_1}^{\ts(r_1)}
\dots
T_{i_mj_m}^{\ts(r_m)}
+X
,
T_{k_1l_1}^{\ts(-s_1)}
\dots
T_{k_nl_n}^{\ts(-s_n)}
+Y
\bigr\rangle
=
\bigl\langle
T_{i_1j_1}^{\ts(r_1)}
\dots
T_{i_mj_m}^{\ts(r_m)}
,
T_{k_1l_1}^{\ts(-s_1)}
\dots
T_{k_nl_n}^{\ts(-s_n)}
\bigr\rangle
$$
for any
$X\in\Y(\gl_N)$ and $Y\in\Y^\ast(\gl_N)$ of degrees respectively
less and more than $s$.

\bpr
\label{prop:L4.2}
For any index $s\geqslant0$, the
restriction of the pairing \eqref{3.93} to\/
${\rm gr}_s\Y(\gl_N)\times{\rm gr}_s\Y^\ast(\gl_N)$ is
non-degenerate.
\epr

\bpf
Fix an integer $s\geqslant0$. In each of
two vector spaces ${\rm gr}_s\Y(\gl_N)$ and ${\rm gr}_s\Y^\ast(\gl_N)$
we will choose a basis so that the matrix of the
bilinear pairing \eqref{3.93}
relative to these bases is lower triangular, with non-zero
diagonal entries.
In particular, we will prove that these two vector spaces are of the
same dimension.

Let $r_1,\dots,r_m$ and $s_1,\dots,s_n$ be
non-increasing sequences of positive integers
satisfying \eqref{self-cons}.
In other words, these two sequences are partitions of
$s$. We will equip the set of all partitions of $s$
with the \textit{inverse} lexicographical ordering. In this
ordering, the sequence $r_1,\dots,r_m$ precedes the sequence
$s_1,\dots,s_n$ if for some $c\geqslant0$
\ben
r_m=s_n\,,\tss
r_{m-1}=s_{n-1}\,,\,
\dots\,,\tss
r_{m-c+1}=s_{n-c+1}
\quad
\text{while}
\quad
r_{m-c}<s_{n-c}\,.
\een

Suppose that $s_1,\dots,s_d\geqslant2$
while $s_{d+1},\dots,s_n=1$ for $d\geqslant 0$.
Unlike in the proof of Lemma~\ref{lem:L4.1}, now
we do not exclude the case $d=n$.
Take the coefficient~at
\beql{3.94}
u_1^{-r_1}\dots u_m^{-r_m}\,v_1^{\ts s_1-1}\dots v_n^{\ts s_n-1}
\eeq
in the expansion of the product \eqref{3.925} as a series in
$u_1^{-1},\dots,u_m^{-1},v_1\,\dots,v_n\tss$.
This coefficient is an element of the algebra
$(\End\CC^N)^{\ot(m+n)}$. If this coefficient is non-zero,
then equality
\ben
r_1+\dots+r_m=s_1+\ldots+s_n
\een
implies that each of the indices $r_1,\dots,r_m$ in \eqref{3.94}
is a sum of some of the indices $s_1,\dots,s_n\ts$. Moreover,
then each of the indices $s_1,\dots,s_n$ appears in these sums only
once. If a sequence $r_1,\dots,r_m$ obtained by this summation
precedes the sequence $s_1,\dots,s_n$ in the inverse lexicographical
ordering, then
the two sequences must coincide. That is,
$m=n$ and $r_a=s_a$ for every index
$a=1,\dots,m$.

For $r=1,2,\dots$ denote by $\Sc_r$ the segment of the sequence
$1,\dots,m$ consisting of all indices $a$ such that $s_a=r$.
If the sequences $r_1,\dots,r_m$ and $s_1,\dots,s_n$ coincide,
then the coefficient at \eqref{3.94} in the expansion of
the product \eqref{3.925} equals
\ben
(-1)^m\ \prod_{r\geqslant1}
\,\,\Bigl(\,\,
\sum_p
\prod_{a\in\Sc_r}
P_{\ts a,p(a)+m}
\Bigr)
\een
where $p$ runs through the set of all permutations of the sequence
$\Sc_r\tss$.
Note that in the products over
$r$ and $a$ above, all the factors pairwise commute.

The graded algebra $\gr\Y(\gl_N)$
is free commutative with the generators
$\Tb_{ij}^{\ts(r)}$ where $r\geqslant1\,$,
see Corollary~\ref{cor:polalg}.
Choose the basis in the vector space
${\rm gr}_s\Y(\gl_N)$ consisting
of the monomials
\beql{grb}
\Tb_{i_1j_1}^{\ts(r_1)}
\dots\,
\Tb_{i_mj_m}^{\ts(r_m)}.
\eeq
The ordering of factors in \eqref{grb} is irrelevant,
let us order them in any way such that
$r_1\geqslant\dots\geqslant r_m$.
Choose any linear ordering of these basis vectors,
subordinate to the inverse lexicographical ordering of the
corresponding sequences $r_1,\dots,r_m\tss$.
The above arguments imply,
that for any two basis elements,
\ben
\Tb_{i_1j_1}^{\ts(r_1)}
\dots\,
\Tb_{i_mj_m}^{\ts(r_m)}
\fand
\Tb_{k_1l_1}^{\ts(s_1)}
\dots\,
\Tb_{k_ml_m}^{\ts(s_m)}
\een
such that the the sequence $r_1,\dots,r_m$ precedes the sequence
$s_1,\dots,s_n$,
the pairing value \eqref{gpvalue}
is non-zero only if $m=n$ and
for every index $a=1,\dots,m$ we have
\ben
r_a=s_a
\fand
i_a=l_a\,,\,
j_a=k_a\,.
\een
Then the value \eqref{gpvalue} equals $(-1)^m\tss g!\tss h!\tss\dots$
where $g,h,\dots$ are the multiplicities in the sequence
of the triples
$
(r_1,i_1,j_1),\dots,(r_m,i_m,j_m).
$
Therefore the monomials
\beql{grbast}
\Tt_{j_1i_1}^{\ts(-r_1)}
\dots\,
\Tt_{j_mi_m}^{\ts(-r_m)}
\eeq
in ${\rm gr}_s\Y^\ast(\gl_N)$
corresponding to the
basis elements \eqref{grb} of vector space ${\rm gr}_s\Y(\gl_N)\,$,
are linearly independent. These monomials also span
the vector space ${\rm gr}_s\Y^\ast(\gl_N)\,$. The latter result
follows from Lemma~\ref{lem:grel} by using induction on $m\,$.
Hence these monomials 
form a basis in ${\rm gr}_s\Y^\ast(\gl_N)\,$.
The matrix of the pairing \eqref{3.93} relative to the two bases
is then lower triangular, with non-zero diagonal entries.
\epf

The graded algebra $\gr\Y^\ast(\gl_N)$
inherits from $\Y^\ast(\gl_N)$ the bialgebra structure. Namely,
using \eqref{Deltamadual}, for any $r\geqslant1$ we get
\beql{gradeddualhopf}
\Delta\bigl(\,\Tt_{ij}^{\ts(-r)}\bigr)=
\Tt_{ij}^{\ts(-r)}\ot1+1\ot\Tt_{ij}^{\ts(-r)}
\fand
\ve\bigl(\,\Tt_{ij}^{\ts(-r)}\bigr)=0.
\eeq
Although the antipode $\Sr$ is defined only
on the completion $\Y^\circ(\gl_N)$ of $\Y^\ast(\gl_N)$,
it still induces a well-defined anipodal map
on the graded algebra $\gr\Y^\ast(\gl_N)$,
\beql{gradedduals}
\Sr:\Tt_{ij}^{\ts(-r)}\mapsto-\,\Tt_{ij}^{\ts(-r)}.
\eeq
Hence $\gr\Y^\ast(\gl_N)$ becomes a Hopf algebra.

Now consider the subalgebra $z\,\gl_N[z]\cong\gl_N\ot(\tss z\,\CC[z]\tss)$
in the polynomial current Lie algebra $\gl_N[z]$. 
The next proposition
indicates the difference between
the graded algebras $\gr\Y(\gl_N)$ and $\gr\Y^\ast(\gl_N)$,
cf.\ Proposition~\ref{prop:iso2}.

\bpr
\label{prop:C4.3}
The Hopf algebra $\gr\Y^\ast(\gl_N)$ is isomorphic to the
universal enveloping algebra 
$\U(\tss z\,\gl_N[z]\tss)$.
\epr

\bpf
Lemma~\ref{lem:grel} implies that the assignment
$E_{ij}\ts z^r\mapsto\Tt_{ij}^{\ts(-r)}$ for $r\geqslant1$
defines a surjective homomorphism
\beql{starhomo}
\U(\tss z\,\gl_N[z]\tss)\to\gr\Y^\ast(\gl_N)\tss.
\eeq
The kernel of this homomorphism is trivial,
because the monomials \eqref{grbast} in
$\Tt_{ij}^{\ts(-r)}$ corresponding to
basis elements \eqref{grb} of the free commutative algebra
$\gr\Y(\gl_N)$ form a basis in $\gr\Y^\ast(\gl_N)\,$.
This was shown in the proof of Proposition~\ref{prop:L4.2}.
By comparing the definitions \eqref{gradeddualhopf},\eqref{gradedduals}
with 
\eqref{standardhopf1},\eqref{standardhopf2} we complete the proof.
\epf

We state the main property of the pairing
$\langle\,\ts,\ts\rangle$ as the following theorem.

\bth
\label{thm:P4.4}
The map \eqref{pairing}
is a non-degenerate bialgebra pairing.
\eth

\bpf
By Lemma~\ref{lem:L4.1} and Proposition~\ref{prop:L4.2}
the pairing $\langle\,\ts,\ts\rangle$
is non-degenerate.
Let us show that under the pairing \eqref{pairing},
the multiplication and comultiplication on $\Y(\gl_N)$
become dual respectively to the comultiplication and multiplication
on $\Y^\ast(\gl_N)$.
We have to prove that
\beql{3.18}
\langle\ts X,Z\ts W\ts\rangle=
\langle\ts\Delta(X),Z\ot W\ts\rangle
\fand
\langle\ts X\ts Y,Z\ts\rangle=
\langle\ts X\ot Y,\Delta(Z)\tss\rangle
\eeq
for any elements $X,Y\in\Y(\gl_N)$ and $Z,W\in\Y^\ast(\gl_N)$.
Here we use the convention
\ben
\langle\ts X\ot Y,Z\ot W\ts\rangle=
\langle\ts X,Z\rangle\ts
\langle\ts Y,W\ts\rangle\ts.
\een

For instance, let us prove the first equality in \eqref{3.18}.
To this end it suffices to substitute the series
$
T_{i_1j_1}(u_1)\dots T_{i_mj_m}(u_m)
$
and
\ben
T^{\tss\ast}_{k_1l_1}(v_1)\dots T^{\tss\ast}_{k_dl_d}(v_d)
\ ,\
T^{\tss\ast}_{k_{d+1}l_{d+1}}(v_{d+1})\dots T^{\tss\ast}_{k_nl_n}(v_n)
\een
for 
$X$ and $Z,W$ respectively.
Here $0\leqslant d\leqslant n$. If $d=0$ or $d=n$, then
we substitute $1$ respectively for $Z$ or for $W$.
After these substitutions, we will have to prove that
\beql{subleft}
\bigl\langle\,
T_{i_1j_1}(u_1)\dots T_{i_mj_m}(u_m)\,,\,
T_{k_1l_1}^{\tss\ast}(v_1)\dots T_{k_nl_n}^{\tss\ast}(v_n)
\,\bigl\rangle
\eeq
equals the sum
\begin{align}
\label{subright}
\sum_{h_1,\dots,h_m=1}^N
&\bigl\langle\,
T_{i_1h_1}(u_1)\dots T_{i_mh_m}(u_m)\,,\,
T_{k_1l_1}^{\tss\ast}(v_1)\dots T_{k_dl_d}^{\tss\ast}(v_n)
\,\bigl\rangle\,\times
\\
\nonumber
&\bigl\langle\,
T_{h_1j_1}(u_1)\dots T_{h_mj_m}(u_m)\,,\,
T_{k_{d+1}l_{d+1}}^{\tss\ast}(v_1)\dots T_{k_nl_n}^{\tss\ast}(v_n)
\,\bigl\rangle\,.
\end{align}

To prove the latter equality,
let us multiply \eqref{subleft} and \eqref{subright} by
the element
\ben
e_{i_1j_1}\ot\dots\ot e_{i_mj_m}
\ot
e_{k_1l_1}\ot\dots\ot e_{k_nl_n}
\in
(\End\CC^N)^{\ot(m+n)},
\een
taking the sum over the indices
$i_1,j_1,\dots,i_m,j_m$ and $k_1,l_1,\dots,k_n,l_n\tss$.
In this way,
from \eqref{subleft} we obtain the product
\ben
\prod_{1\leqslant a\leqslant m}^\rightarrow
\,
\prod_{1\leqslant b\leqslant n}^\rightarrow R_{a,b+m}(u_a-v_b)
\een
due to the definition \eqref{3.7777777}.
From \eqref{subright} we obtain the product
\ben
\prod_{1\leqslant b\leqslant d}^\rightarrow
\,
\prod_{1\leqslant a\leqslant m}^\rightarrow
R_{a,b+m}(u_a-v_b)
\ \cdot\!
\prod_{d<b\leqslant n}^\rightarrow
\,
\prod_{1\leqslant a\leqslant m}^\rightarrow
R_{a,b+m}(u_a-v_b)
\een
which is evidently equal to the previous product.

We have already noted the equality $\langle\ts1\ts,1\,\rangle=1$.
Moreover, by setting $n=0$
the definition \eqref{3.7777777}, for any
$r_1,\dots,r_m\geqslant1$ we get the equality
\ben
\langle\ts
T_{i_1j_1}^{(r_1)}\dots T_{i_mj_m}^{(r_m)}\tss,\tss1
\ts\rangle=0
\quad\text{if}\quad
m\geqslant1\,.
\een
Thus $\langle\tss X\tss,1\,\rangle=\ve(X)$
for any element $X\in\Y(\gl_N)$.
By setting $m=0$ in
\eqref{3.7777777} and using the induction on $n$
or, alternatively, by using Lemma~\ref{lem:L4.1},
we obtain for any $s_1,\dots,s_n\geqslant1$ the equality
\ben
\langle\ts1\ts,\tss
T_{k_1l_1}^{\ts(-s_1)}\dots T_{k_nl_n}^{\ts(-s_n)}
\ts\rangle=0
\quad\text{if}\quad
n\geqslant1.
\een
Thus $\langle\ts1\ts,Z\,\rangle=\ve(Z)$
for any element $Z\in\Y^\ast(\gl_N)\,$.
Therefore the counit and the unit maps for the bialgebra
$\Y(\gl_N)$ are dual respectively to the unit and the counit
maps for the bialgebra $\Y^\ast(\gl_N)$.
\epf

Due to Theorem~\ref{thm:pbw},
the vector space $\Y(\gl_N)$ has a basis parameterized
by all multisets of triples $(r_1,i_1,j_1),\dots,(r_m,i_m,j_m)$ where
\ben
r_1,\dots,r_m\in\{1,2,\dots\}
\fand
i_1,j_1,\dots,i_m,j_m\in\{1,\dots,N\}
\een
while $m=0,1,2,\dots\ts\tss$.
The corresponding basis vector in $\Y(\gl_N)$ is the monomial
\beql{basis}
T_{i_1j_1}^{\ts(r_1)}
\dots\,
T_{i_mj_m}^{\ts(r_m)}.
\eeq
The ordering of the factors in this monomial can be chosen
arbitrarily. Suppose that here
$r_1\geqslant\dots\geqslant r_m\tss$.
Then the sequence $r_1,\dots,r_m$ can be regarded as a partition of
$r_1+\dots+r_m\tss$.
Equip the set of all partitions of $0,1,2,\dots$
with the following ordering. If $r<s$, the partitions of
$r$ precede those of $s$. For any given $r$,
the set of partitions of $r$ is equipped with the
inverse lexicographical ordering\tss;
see the proof of Proposition~\ref{prop:L4.2}. Choose any linear
ordering of the basis elements \eqref{basis},
subordinate to the above described ordering of
their sequences $r_1,\dots,r_m\tss$.
The proof of Proposition~\ref{prop:L4.2} implies that the monomials
\ben
T_{j_1i_1}^{\ts(-r_1)}
\dots\,
T_{j_mi_m}^{\ts(-r_m)}
\een
corresponding to the basis elements \eqref{basis}
form a basis of the vector space $\Y^\ast(\gl_N)$.
The matrix of the 
pairing \eqref{pairing}
relative to these two bases is lower triangular with non-zero
diagonal entries; see also Lemma~\ref{lem:L4.1}.
Here the basis elements of $\Y^\ast(\gl_N)$
are linearly ordered as the corresponding basis elements
\eqref{basis} of $\Y(\gl_N)\,$.


\section{Universal $R$-matrix}
\label{sec:unir}

Consider the formal completion $\Y^\circ(\gl_N)$ of the
filtered algebra $\Y^\ast(\gl_N)$ defined in Section~\ref{sec:dualyang}.
By Proposition \ref{lem:L4.1}
the canonical pairing \eqref{pairing} extends to a pairing
\ben
\langle\,\ts,\,\rangle:\Y(\gl_N)\times\Y^\circ(\gl_N)\to\CC\tss.
\een
Choose any basis $X_1,X_2,\dots$ in the vector space $\Y(\gl_N)$.

\bpr
\label{prop:dualsys}
The completion $\Y^\circ(\gl_N)$ does
contain the system of elements
$X'_1,X'_2,\dots$ dual to $X_1,X_2,\dots$ so that
$\langle\,X_r\tss,\tss X'_s\tss\rangle=\de_{rs}$
for any $r$ and $s$.
\epr

\bpf
As we explained at the end of Section \ref{sec:dualyang},
one can choose a basis $Y_1,Y_2,\dots$ in $\Y(\gl_N)$ and
a basis $Y^\ast_1,Y^\ast_2,\dots$ in $\Y^\ast(\gl_N)$ so that
the filtration degree
\beql{yasttoinfty}
\deg\,Y_s^\ast\to\infty
\quad\text{when}\quad
s\to\infty\tss,
\eeq
and so that the matrix
of the pairing \eqref{pairing} relative
to these bases is lower triangular with non-zero diagonal entries.
Let $[\,g_{rs}\tss]$ be its inverse matrix. The formal sums
\beql{gsum}
Y'_s=\sum_{r=1}^\infty\,g_{\tss rs}\,Y^\ast_r
\eeq
satisfy the equations
$\langle\,Y_r\tss,\tss Y'_s\tss\rangle=\de_{rs}$
for all indices $r$ and $s$.
Each of these sums is contained in $\Y^\circ(\gl_N)$
due to \eqref{yasttoinfty}.
Moreover, because the the matrix $[g_{rs}]$ is also lower triangular,
the property \eqref{yasttoinfty} implies that
\beql{yprimetoinfty}
\deg Y_s^\prime\to\infty
\quad\text{when}\quad
s\to\infty\tss.
\eeq

Now let $X_1,X_2,\dots$ be any basis in $\Y(\gl_N)$. Let
$[h_{rs}]$ be the coordinate change matrix from the basis
$Y_1,Y_2,\dots$ so that for any index $r$ we have
\ben
Y_s=\sum_{r=1}^\infty\,h_{\tss rs}\,X_r\,.
\een
This sum must be finite, so that for any fixed index $s$
there are only finitely many non-zero coefficients $h_{rs}\tss$.
The sums
\beql{xprime}
X'_r=\sum_{s=1}^\infty\,h_{\tss rs}\,Y'_s
\eeq
satisfy the equations
$\langle\,X_r\tss,\tss X'_s\tss\rangle=\de_{rs}$ as required.
Each of these sums is contained in the completion $\Y^\circ(\gl_N)$
due to the property \eqref{yprimetoinfty}.
\epf

Consider an infinite sum of elements of the tensor product
$\Y^\circ(\gl_N)\ot\Y(\gl_N)$
\beql{unir}
\Rc=\sum_{r=1}^\infty\,X'_r\ot X_r\,.
\eeq
This sum does not depend
on the choice of the basis $X_1,X_2,\dots$ in the vector space $\Y(\gl_N)$
in the following sense. Let $Y_1,Y_2,\dots$ be the basis in $\Y(\gl_N)$
used in the proof of Proposition~\ref{prop:dualsys}.
Using the formula \eqref{xprime} for every $r=1,2,\dots$
expand the vectors $X'_1,X'_2,\dots$ in \eqref{unir}.
Then fix an index $s$ and consider the sum of terms
\ben
\sum_{r=1}^\infty\,
(\tss h_{\tss rs}\,Y'_s\tss)\ot X_r
\een
corresponding to the vector $Y'_s$ in \eqref{xprime}.
Only finite number of these terms are non-zero,
and their sum is equal to $Y'_s\ot Y_s\tss$. In this
sense, the 
sum in \eqref{unir} equals
\ben
\sum_{s=1}^\infty\,
Y'_s\ot Y_s\,.
\een
The infinite sum $\Rc$ 
is called the {\it universal R-matrix} for the Yangian $\Y(\gl_N)\,$.

Any element of the vector space $\Y^\circ(\gl_N)\ot\Y(\gl_N)$
determines a linear operator on the vector space $\Y(\gl_N)\,$. If $A$
is the operator corresponding to an element
$Z\ot Y\in\Y^\circ(\gl_N)\ot\Y(\gl_N)$, then
\beql{pointwise}
A\tss(X)=\langle\ts X\tss,Z\ts\rangle\tss\,Y
\quad\ \text{for any}\quad
X\in\Y(\gl_N).
\eeq
By the above argument,
the series of operators corresponding to \eqref{unir}
pointwise converges to the identity operator
${\id}:X\mapsto X$ on the vector space $\Y(\gl_N)\,$.

\bpr
\label{prop:comr}
For the comultiplication on $\Y(\gl_N)$ and\/ $\Y^\circ(\gl_N)$ we have
\beql{comr}
(\tss{\id}\ot\Delta)\,(\Rc)=\Rc_{12}\tss\Rc_{13}
\fand\/
(\tss\Delta\ot{\id})\,(\Rc)=\Rc_{13}\tss\Rc_{23}
\eeq
where
\ben
\Rc_{12}=\sum_{r=1}^\infty\,X'_r\ot X_r\ot 1\,,
\quad
\Rc_{13}=\sum_{r=1}^\infty\,X'_r\ot 1\ot X_r\,,
\quad
\Rc_{23}=\sum_{r=1}^\infty\,1\ot X'_r\ot X_r\,.
\een
\epr

\bpf
Let us prove the first of the two equalities \eqref{comr}.
This is an equality of infinite sums of elements from the
tensor product $\Y^\circ(\gl_N)\ot\Y(\gl_N)\ot\Y(\gl_N)\tss$.
It means the equality of the corresponding
operators $\Y(\gl_N)\to\Y(\gl_N)\ot\Y(\gl_N)\tss$.
By applying the linear operator corresponding to the infinite sum
$(\tss{\id}\ot\Delta)\,(\Rc)$ to any fixed element $X\in\Y(\gl_N)$
we get the element $\Delta(X)$.
By applying to $X$ the operator corresponding to $\Rc_{12}\tss\Rc_{13}$
we obtain the sum
\ben
\sum_{r,s=1}^\infty\,
\langle\,X\,,Y'_r\,\tss Y'_s\,\,\rangle\,Y_r\ot Y_s
=
\sum_{r,s=1}^\infty\,
\langle\,\Delta(X)\,,\tss Y'_r\ot Y'_s\,\,\rangle\,Y_r\ot Y_s
=
\Delta(X).
\een
Here we used the first equality in \eqref{3.18},
and non-degeneracy of the pairing \eqref{pairing}.
The property \eqref{yprimetoinfty} guarantees that
in both sums over $r$ and $s$ displayed above, only finite number
of summands are non-zero when $X$ is fixed; see Lemma~\ref{lem:L4.1}.
We have thus proved the first equality in \eqref{comr}.
The second equality
is deduced from the second equality in \eqref{3.18}
in a similar way.
\epf

\bpr
\label{prop:cour}
For the counit maps on $\Y(\gl_N)$ and $\Y^\circ(\gl_N)$,
\ben
(\tss{\id}\ot\ve)\ts(\Rc)=1
\fand
(\tss\ve\ot{\id})\ts(\Rc)=1.
\een
\epr

\bpf
Because $\ve(X)=\langle\,X\,,1\,\rangle$
for any element $X\in\Y(\gl_N)$ by Theorem~\ref{thm:P4.4},
\ben
(\tss{\id}\ot\ve)\ts(\Rc)=
\sum_{s=1}^\infty\,
\langle\,Y_s\,,1\,\rangle\,Y'_s
=1\tss.
\een
Similarly, because $\ve(Z)=\langle\,1\,,Z\,\rangle$
for any element $Z\in\Y^\circ(\gl_N)$, we also have
\ben
(\tss\ve\ot{\id})\ts(\Rc)=
\sum_{s=1}^\infty\,
\langle\,1\,,\tss Y'_s\,\rangle\,Y_s=1
\een
where only finitely
many summands are non-zero due to 
\eqref{yprimetoinfty},
see Lemma~\ref{lem:L4.1}.
\epf

The infinite sum in \eqref{unir} can be also regarded
as an element of a 
completion of the tensor product $\Y^\ast(\gl_N)\ot\Y(\gl_N)$. Namely,
let us extend the descending filtration from the algebra $\Y^\ast(\gl_N)$
to the tensor product $\Y^\ast(\gl_N)\ot\Y(\gl_N)$ by giving the
degree $r$ to each element of the form $T_{ij}^{\ts(-r)}\ot X$
where $X\in\Y(\gl_N)$ and $r\geqslant1$. The element
$1\ot X$ of $\Y^\ast(\gl_N)\ot\Y(\gl_N)$ is given the zero degree. Take
the formal completion of the algebra $\Y^\ast(\gl_N)\ot\Y(\gl_N)$
relative to this filtration.
This completion contains the tensor product
$\Y^\circ(\gl_N)\ot\Y(\gl_N)\tss$, but does not coincide with it
because the algebra $\Y(\gl_N)$ is infinite-dimensional.

The next corollary shows in particular, that the sum in \eqref{unir}
is invertible as an element of the completion of the algebra
$\Y^\ast(\gl_N)\ot\Y(\gl_N)\tss$.

\bco
\label{cor:antir}
For the antipodal maps on $\Y(\gl_N)$ and\/ $\Y^\circ(\gl_N)$ we have
\ben
(\tss{\id}\ot\Sr)\ts(\Rc)=\Rc^{-1}
\fand\/
(\tss\Sr\ot{\id})\ts(\Rc)=\Rc^{-1}.
\een
\eco

\bpf
Regard the
first equality in \eqref{comr} as that of the elements
of the completion of the algebra
$\Y^\ast(\gl_N)\ot\Y(\gl_N)\ot\Y(\gl_N)$.
On this algebra, the
descending filtration is defined
by giving the degree $r$ to each element of the form
$T_{ij}^{\ts(-r)}\ot X\ot Y$ where
$X,Y\in\Y(\gl_N)$ and $r\geqslant1$.
The element $1\ot X\ot Y$ is then given the degree zero.

Let $\mu:\Y(\gl_N)\ot\Y(\gl_N)\to\Y(\gl_N)$ be the map of algebra
multiplication, and $\de:\CC\to\Y(\gl_N)$ be the unit map\tss:
$\de\tss(1)=1\tss$.
Let us apply the map ${\id}\ot\Sr\ot{\id}$, and then
the map ${\id}\ot\mu$ to
to both sides of the first equality in \eqref{comr}.
At the right hand side we get the element
$((\tss{\id}\ot\Sr)\tss(\Rc))\cdot\Rc\tss$.
At the left hand side we get the element of the tensor product
$\Y^\circ(\gl_N)\ot\Y(\gl_N)\tss$,
\ben
(\tss(\tss{\id}\ot\mu)\,
(\tss{\id}\ot\Sr\ot{\id})\,
(\tss{\id}\ot\Delta)\tss)\,(\Rc)
=
(\tss(\tss{\id}\ot\de)\,
(\tss{\id}\ot\ve)\tss)\,
(\Rc)=1\ot1\tss.
\een
Here we used the first axiom of antipode
from Section~\ref{sec:hopf}
in the case $\Ar=\Y(\gl_N)$,
and the first equality of Proposition~\ref{prop:cour}.
Hence the first equality of Corollary~\ref{cor:antir}
follows from the first equality in \eqref{comr}.

Similarly, using the first axiom of antipode
in the case $\Ar=\Y^\circ(\gl_N)$ and
the second equality of Proposition~\ref{prop:cour},
the second equality of Corollary~\ref{cor:antir}
follows from the second equality in \eqref{comr}.
The last equality should be regarded here
as that of the elements of the completion of the algebra
$
\Y^\ast(\gl_N)\ot\Y^\ast(\gl_N)\ot\Y(\gl_N)\tss.
$
On this algebra a
descending filtration is defined
by giving the degree $r+s$ to any element of the form
$$
T_{ij}^{\ts(-r)}\ot T_{kl}^{\ts(-s)}\ot X\,.
$$
Then the elements 
$T_{ij}^{\ts(-r)}\ot1\ot X$ and $1\ot T_{ij}^{\ts(-r)}\ot X$
are given the degree $r$, while
the element $1\ot 1\ot X$ is given the degree zero.
Here $X\in\Y(\gl_N)$ and $r,s\geqslant1\,$.
The argument is completed in the same way as for the first
equality.
\epf

Let us now replace the complex parameter $c$ in the definition \eqref{3.66}
of a covector representation $\rho_c$ of $\Y(\gl_N)$ by the formal
variable $v\tss$. 
Then we get a homomorphism
\beql{rhoc}
\rho_v:
\Y(\gl_N)\to\End\CC^N\tss[v]\tss;
\eeq
it is defined by the assignment $T(u)\mapsto R(u-v)$
of formal power series in $u^{-1}$.

Similarly, the assignment $T^{\tss\ast}(v)\mapsto R(u-v)$
of formal power series in $v$ defines a homomorphism
\beql{rhocast}
\rho^{\tss\ast}_u:
\Y^\ast(\gl_N)\to\End\CC^N\tss[u^{-1}]\tss.
\eeq
To prove the homomorphism property
using the matrix form \eqref{3.13} of the defining relations
of the algebra $\Y^\ast(\gl_N)\tss$,
we have to check the equality of rational functions in
the variables $u,v$ and $w$
with the values in the algebra $(\End \CC^N)^{\ot3}$,
\ben
R_{\tss01}(u-v)\ts R_{\tss02}(u-w)\ts R_{12}(v-w)=
R_{12}(v-w)\ts R_{\tss02}(u-w)\ts R_{\tss01}(u-v)\tss.
\een
This equality follows from \eqref{ybe}.
Here we use the indices $0,1,2$ instead of $1,2,3$
to label the tensor factors of $(\End \CC^N)^{\ot3}$.
By comparing
the expansions \eqref{tutensast} and \eqref{rexp}, we see that
\beql{3.666}
\rho_u^{\tss\ast}:
T_{ij}^{\tss(-r)}\mapsto-\tss u^{\tss-r}e_{ji}
\,\quad\text{for any}\quad
r\geqslant1\tss.
\eeq
Obviously, the homomorphism $\rho_u^{\tss\ast}$
extends to a homomorphism
\ben
\Y^\circ(\gl_N)\to\End\CC^N\tss[[u^{-1}]]\tss.
\een
We shall keep the notation $\rho_u^{\tss\ast}$ for the extended homomorphism.

\bpr
\label{prop:P4.5}
We have equalities of formal power series in $u^{-1}$ and $v$,
\ben
(\tss\rho_u^{\tss\ast}\ot{\id}\tss)\,(\Rc)=T(u)
\fand
(\tss\id\ot\rho_v\tss)\,(\Rc)=T^{\tss\ast}(v)\,.
\een
\epr

\bpf
By the definition \eqref{3.7777777}
of the pairing $\Y(\gl_N)\ot\Y^\ast(\gl_N)\to\CC$ for any $n\geqslant0$
the element $T(u)\in\End(\CC^N)\ot\Y(\gl_N)[[u^{-1}]]$
has the property that
\ben
T(u)\ot T_1^{\tss\ast}(v_1)\dots T_n^{\tss\ast}(v_n)\mapsto
R_{12}(u-v_1)\dots R_{1,n+1}(u-v_n)
\een
under the linear map
\ben
\id\ot\be\ot\id\tss:
\End\CC^N\ot\Y(\gl_N)\ot\Y^\ast(\gl_N)\ot(\End\CC^N)^{\ot n}
\to(\End\CC^N)^{\ot\tss(n+1)}.
\een
Because our pairing is non-degenerate, the same property for
the element
\ben
(\tss\rho_u^\ast\ot{\id}\tss)\,(\Rc)=\sum_{s=1}^\infty\,
\rho_u^\ast\tss(\tss Y'_s)\ot Y_s
\een
will imply the first equality of Proposition \ref{prop:P4.5}.
Note that when $s\to\infty\tss$, then
the degree in $u^{-1}$ of the image
$\rho_u^\ast\tss(\tss Y'_s)$ tends to infinity
due to \eqref{yprimetoinfty} and \eqref{3.666}. Hence
the above displayed sum over $s=1,2,\dots$ is contained in
$\End(\CC^N)\ot\Y(\gl_N)[[u^{-1}]]\,$.

Thus to prove the first equality of Proposition \ref{prop:P4.5},
we have to show that under the linear map $\id\ot\be\ot\id\tss$,
\ben
\sum_{s=1}^\infty\,
\rho_u^{\tss\ast}\tss(\tss Y'_s)\ot Y_s
\ot T_1^{\tss\ast}(v_1)\dots T_n^{\tss\ast}(v_n)\mapsto
R_{12}(u-v_1)\dots R_{1,n+1}(u-v_n)\tss.
\een
Since the system of vectors $Y'_1,Y'_2,\dots$
is dual to the basis $Y_1,Y_2,\dots$ of $\Y(\gl_N)\tss$,
this is equivalent to showing that
\ben
T_1^{\tss\ast}(v_1)\dots T_n^{\tss\ast}(v_n)\mapsto
R_{12}(u-v_1)\dots R_{1,n+1}(u-v_n)
\een
under the linear map
\ben
\rho_u^{\tss\ast}\ot\id:
\Y^\ast(\gl_N)\ot(\End\CC^N)^{\ot n}
\to(\End\CC^N)^{\ot\tss(n+1)}[u^{-1}]\tss.
\een
The latter property follows directly
from the definition of the homomorphism $\rho^{\tss\ast}_u\tss$.
The proof of the second equality of Proposition \ref{prop:P4.5}
is similar and is omitted.
\epf

\bco
\label{cor:r}
We have the equality of formal power series in $u^{-1}$ and $v$,
\ben
(\tss\rho_u^{\tss\ast}\ot\rho_v\tss)\,(\Rc)=R(u-v)\tss.
\een
\eco


\section{Double Yangian}
\label{sec:doubleyang}

Let $\Delta'$ be the comultiplication on $\Y^\ast(\gl_N)$
\textit{opposite} to the comultiplication $\Delta$ defined by
\eqref{3.17}. By definition, the map
$$
\Delta':\Y^\ast(\gl_N)\to\Y^\ast(\gl_N)\ot\Y^\ast(\gl_N)
$$
is the composition of the comultiplication $\Delta$
with the linear operator on
the tensor product $\Y^\ast(\gl_N)\ot\Y^\ast(\gl_N)$ exchanging
the tensor factors.

The \textit{double Yangian\/} of
$\gl_N$ is defined as an associative unital algebra
$\DY(\gl_N)$ over $\CC$
generated by the elements of $\Y(\gl_N)$ and $\Y^\ast(\gl_N)$
subject to the relations
\beql{3.99}
\Rc\,\tss\Delta(W)=\Delta'(W)\,\Rc
\quad\,\text{for every}\quad
W\in\Y^\ast(\gl_N)\,.
\eeq
In the rest of this section we will provide a more explicit
description of the algebra $\DY(\gl_N)\,$, see Theorem~\ref{cor:doublerel}
below. In Section \ref{sec:dpbw} we will show that the defining
homomorphisms of 
$\Y(\gl_N)$ and $\Y^\ast(\gl_N)$ to $\DY(\gl_N)$ are in fact embeddings.
At the end of that section we will also provide an equivalent
definition of the $\DY(\gl_N)\,$.

In \eqref{3.99} we have an equality of infinite sums of
elements of the tensor product $\Y^\circ(\gl_N)\ot\DY(\gl_N)$.
It means the equality of the corresponding linear
operators $\Y(\gl_N)\to\DY(\gl_N)$, cf. \eqref{pointwise}.
For instance, let us consider the infinite sum
\ben
\Rc\,\tss\Delta(W)=\sum_{s=1}^\infty\,(\tss Y'_s\ot Y_s)\,\tss\Delta(W)
\een
at the right hand side of the equality postulated in \eqref{3.99}.
Note that for any fixed
$X\in\Y(\gl_N)$ and $Z\in\Y^\ast(\gl_N)$,
only finitely many summands in the infinite sum
\ben
\sum_{s=1}^\infty\,
\langle\,X\,,Y'_s\tss Z\,\rangle\,Y_s
\een
are non-zero; see Lemma~\ref{lem:L4.1} and the property
\eqref{yprimetoinfty}. This observation shows that
the linear operator $\Y(\gl_N)\to\DY(\gl_N)$ corresponding
to the infinite sum $\Rc\tss\,\Delta(W)$ is well-defined for any
element $W\in\Y^\ast(\gl_N)\tss$.

Now take the pair of homomorphisms $\rho_u$ and $\rho^{\tss\ast}_u$
where we use the same formal variable $u\tss$,
see \eqref{rhoc} and \eqref{rhocast}.

\bpr
\label{prop:C4.6}
The associative algebra homomorphisms\/ $\rho_u\tss,\rho^{\tss\ast}_u$
extend to a homomorphism
$
\DY(\gl_N)\to\End\CC^N[\ts u,u^{-1}]\,.
$
\epr

\bpf
Using \eqref{3.99}, for any $W\in\Y^\ast(\gl_N)$
we have to check the equality
\ben
(\tss\id\ot\rho_u\tss)\,(\Rc)
\,
(\tss\id\ot\rho_u^\ast\tss)\,(\Delta\tss(W))
=
(\tss\id\ot\rho_u^\ast\tss)\,(\Delta'(W))
\,
(\tss\id\ot\rho_u\tss)\,(\Rc)
\een
of formal series in $u$ with coefficients in
the algebra $\Y^\circ(\gl_N)\ot\End\CC^N$.
It suffices to substitute here the series $T_{ij}^{\tss\ast}(v)$
for the element $W$.
Due to the definition \eqref{3.17}
and to Proposition \ref{prop:P4.5},
the result of the substitution is the relation
\ben
\sum_{k=1}^N\,\,
T^{\tss\ast}(u)
\,
(\tss T^{\tss\ast}_{ik}(v)\ot\rho_u^\ast(\tss T^{\tss\ast}_{kj}(v)))
=
\sum_{k=1}^N\,\,
(\tss T^{\tss\ast}_{kj}(v)\ot\rho_u^\ast(\tss T^{\tss\ast}_{ik}(v)))
\,
T^{\tss\ast}(u)\tss.
\een
Let us take the tensor products of both sides of the latter relation
with the element $e_{ij}\in\End\CC^N$, and then sum over
$i,j=1\dots,N$. Using the identity
$e_{ij}=e_{ik}\tss e_{kj}$ we then get the relation
\begin{align}
\label{lastrel}
&\sum_{i,j,k=1}^N
(\tss T^{\tss\ast}(u)\ot1\tss)
\,
(\tss T^{\tss\ast}_{ik}(v)\ot\rho_u^\ast(\tss T^{\tss\ast}_{kj}(v))
\ot e_{ik}\tss e_{kj}\tss)\ =
\\
\nonumber
&\sum_{i,j,k=1}^N
(\tss T^{\tss\ast}_{kj}(v)\ot\rho_u^\ast(\tss T^{\tss\ast}_{ik}(v))
\ot e_{ik}\tss e_{kj}\tss)
\,
(\tss T^{\tss\ast}(u)\ot1\tss)
\end{align}
of formal power series in $u,v$ with the coefficients
in $\Y^\ast(\gl_N)\ot\End\CC^N\ot\End\CC^N$.
Note that by  the definition of the homomorphism \eqref{rhocast},
\ben
\sum_{i,j=1}^N
\rho_u^\ast(\tss T^{\tss\ast}_{ij}(v))\ot e_{ij}=R(u-v)\tss.
\een
Therefore the relation \eqref{lastrel} can be rewritten as
\ben
T^{\tss\ast}_1(u)\ts T^{\tss\ast}_2(v)\,(\tss1\ot R(u-v))=
(\tss1\ot R(u-v))\,T^{\tss\ast}_2(v)\ts T^{\tss\ast}_1(u)\tss.
\een
But this is just the defining relation for
the algebra $\Y^\ast(\gl_N)\tss$, see \eqref{3.13}.
\epf

Let $c$ be any non-zero complex number.
In Proposition~\ref{prop:C4.6},
we can specialize the formal variable $u$ to $c\tss$. 
Then we obtain a representation $\DY(\gl_N)\to\End\CC^N$. 
We call it a \textit{covector representation\/}
of the algebra $\DY(\gl_N)$, 
it extends the covector representation \eqref{3.66} of
the algebra $\Y(\gl_N)$. 

The vector representation \eqref{3.666666}
of $\Y(\gl_N)$ can be extended to a representation of
$\DY(\gl_N)\tss$, by mapping $T^{\tss\ast}(v)\mapsto R^{\tss t}(v-u)\tss$.
We call it a \textit{vector representation\/}
of the algebra $\DY(\gl_N)$ and denote it by 
$\sigma_c\tss$. Note that then
\beql{dualcovec}
\si_c:
T_{ij}^{(-r)}\mapsto c^{\tss -r}\tss e_{ij}
\,\quad\text{for any}\quad
r\geqslant1\tss.
\eeq
The proof that these assignments together with \eqref{3.666666}
define a representation of the algebra $\DY(\gl_N)$ is similar to
that of Proposition~\ref{prop:C4.6}, and is omitted here.

To write down commutation relations in the algebra $\DY(\gl_N)\tss$,
we will use the tensor product
$
\End\CC^N\ot\DY(\gl_N)\ot\End\CC^N.
$
There is a natural embedding
of the algebra $\End\CC^N\ot\End\CC^N$ into this tensor product,
such that $x\ot y\mapsto x\ot1\ot y$ for any elements $x,y\in\End\CC^N$.
In the next proposition, 
the Yang $R$-matrix \eqref{rmatrix} is identified with its image
relative to this embedding.

\bpr
\label{prop:C4.7}
In the algebra
$\End\CC^N\ot\DY(\gl_N)\ot\End\CC^N\tss[[\tss u^{-1},v\tss]]$
we have 
\beql{3.111111}
(\tss T(u)\ot1\tss)\,R(u-v)\,(\tss1\ot T^{\tss\ast}(v))
=
(\tss1\ot T^{\tss\ast}(v))\,R(u-v)\,(\tss T(u)\ot1\tss)\,.
\eeq
\epr

\bpf
Let us substitute $T^{\tss\ast}_{ij}(v)$ for $W$ in
the equality in \eqref{3.99},
and then apply the homomorphism $\rho_u^{\tss\ast}\ot\id$ 
to the resulting equality.
Due to the definition \eqref{3.17} and to Proposition \ref{prop:P4.5}, we get
an equality of formal power series in $u^{-1}$ and $v$
with the coefficients from $\End\CC^N\ot\DY(\gl_N)\tss$,
\ben
\sum_{k=1}^N\,\,
T(u)
\,
(\tss\rho_u^{\tss\ast}(\tss T^{\tss\ast}_{ik}(v))\ot T^{\tss\ast}_{kj}(v))
=
\sum_{k=1}^N\,\,
(\tss\rho_u^{\tss\ast}(\tss T^{\tss\ast}_{kj}(v))\ot T^{\tss\ast}_{ik}(v))
\,
T(u)\,.
\een
Let us now take the tensor products of both sides of this equality
with the element $e_{ij}\in\End\CC^N$, and then sum over
$i,j=1\dots,N$. Using the identity
$e_{ij}=e_{ik}\tss e_{kj}$ we obtain
an equality of series
with coefficients from $\End\CC^N\ot\DY(\gl_N)\ot\End\CC^N$
\begin{align}
\label{lasteq}
&\sum_{i,j,k=1}^N
(\tss T(u)\ot1\tss)
\,
(\tss\rho_u^{\tss{\tss\ast}}(\tss T^{\tss{\tss\ast}}_{ik}(v))\ot 
T^{\tss\ast}_{kj}(v)
\ot e_{ik}\tss e_{kj}\tss)\ =
\\
\nonumber
&\sum_{i,j,k=1}^N
(\tss\rho_u^{\tss\ast}(\tss T^{\tss\ast}_{kj}(v))\ot T^{\tss\ast}_{ik}(v)
\ot e_{ik}\tss e_{kj}\tss)
\,
(\tss T(u)\ot1\tss)\tss.
\end{align}
By using the definition of $\rho_u^{\tss\ast}$
the equality \eqref{lasteq} can be rewritten as \eqref{3.111111}.
\epf

\bpr
\label{thm:T4.8}
Relation 
\eqref{3.111111} is equivalent to the 
collection of relations \eqref{3.99}.\!
\epr

\bpf
By Proposition \ref{prop:C4.7}
the relation \eqref{3.111111} follows from \eqref{3.99}.
Let $u_1,\dots,u_m$ be independent variables.
Define the homomorphism
\beql{rhoccast}
\rho^{\tss\ast}_{u_1\dots u_m}:
\Y^\ast(\gl_N)
\to
(\End\CC^N)^{\ot m}\tss[u_1^{-1},\dots,u_m^{-1}]
\eeq
as the composition of the $m$-fold comultiplication
$\Y^\ast(\gl_N)\to\Y^\ast(\gl_N)^{\ot m}$ and of~the
tensor product of the homomorphisms \eqref{rhocast}
where $u=u_1,\dots,u_m\tss$.
By using the descending filtration on
$\Y^\ast(\gl_N)$ and the surjective homomorphism
\eqref{starhomo} we can prove that when the number $m$ vary,
the kernels of all homomorphisms
$\rho^{\tss\ast}_{u_1\dots u_m}$
have only zero intersection.
The proof is similar that of Proposition~\ref{thm:P2.2}
and is omitted here. It now suffices to derive from
\eqref{3.111111} that for any $W\in\Y^{\tss\ast}(\gl_N)$
\beql{3.999999999}
(\tss\rho^{\tss\ast}_{u_1\dots u_m}\!\ot\id\tss)\,(\tss\Rc\,\tss\Delta(W))=
(\tss\rho^{\tss\ast}_{u_1\dots u_m}\!\ot\id\tss)\,
(\tss\Delta'(W)\,\tss\Rc\tss)
\tss.
\eeq
Here the homomorphism \eqref{rhoccast} is extended to a homomorphism
\ben
\Y^\circ(\gl_N)
\to
(\End\CC^N)^{\ot m}\tss[[u_1^{-1},\dots,u_m^{-1}]]
\een
and the extension is still denoted by
$\rho^{\tss\ast}_{u_1\dots u_m}$. Using Propositions \ref{prop:comr}
and \ref{prop:P4.5}, the relation \eqref{3.999999999}
can be rewritten as
\begin{gather*}
T_1(u_1)\ldots T_m(u_m)
\,
(\tss\rho^{\tss\ast}_{u_1\dots u_m}\!\ot\id\tss)\,(\tss\Delta(W))
\\=
(\tss\rho^{\tss\ast}_{u_1\dots u_m}\!\ot\id\tss)\,(\tss\Delta'(W))
\,
T_1(u_1)\ldots T_m(u_m)\tss.
\end{gather*}
It suffices to verify the latter relation for each of the series
$T_{ij}^{\tss\ast}(v)$ being substituted for the element $W$.
By the definition \eqref{3.17},
the substitution yields the relation
of the formal power series in $u_1^{-1},\dots,u_m^{-1}$ and $v$
with the coefficients in the algebra $(\End\CC^N)^{\ot m}\ot\DY(\gl_N)\tss$,
\begin{gather*}
T_1(u_1)\ldots T_m(u_m)\ \times
\\
\sum_{k_1,\dots,k_m=1}^N
\rho^{\tss\ast}_{u_1}(\tss T^{\tss\ast}_{ik_1}(v))
\ot
\rho^{\tss\ast}_{u_2}(\tss T^{\tss\ast}_{k_1k_2}(v))
\ot\dots\ot
\rho^\ast_{u_m}(\tss T^{\tss\ast}_{k_{m-1}k_m}(v))
\ot
T_{k_mj}^{\tss\ast}(v)\ =
\\
\sum_{k_1,\dots,k_m=1}^N
\rho^{\tss\ast}_{u_1}(\tss T^{\tss\ast}_{k_1k_2}(v))
\ot\dots\ot
\rho^{\tss\ast}_{u_{m-1}}(\tss T^{\tss\ast}_{k_{m-1}k_m}(v))
\ot
\rho^{\tss\ast}_{u_m}(\tss T^{\tss\ast}_{k_mj}(v))
\ot
T_{ik_1}^\ast(v)
\\[4pt]
\times\ \
T_1(u_1)\ldots T_m(u_m)\tss.
\end{gather*}

\vspace{4pt}\noindent
Let us now take the tensor products of both sides of this relation
with the element $e_{ij}\in\End\CC^N$, and then sum over the indices
$i,j=1\dots,N$. By using the identity
\ben
e_{ij}=e_{ik_1} e_{k_1k_2}\dots e_{k_{m-1}k_m}e_{k_mj}
\een
in $\End\CC^N$,
we arrive at the following relation of series
with the coefficients from the tensor product
$(\End\CC^N)^{\ot m}\ot\DY(\gl_N)\ot\End\CC^N$:
\begin{gather*}
(\tss T_1(u_1)\dots T_m(u_m)\ot1\tss)
\,
R_{1,m+1}(u_1-v)\dots R_{m,m+1}(u_m-v)
\,
(\tss1\ot T^{\tss\ast}(v))
\\=\,
(\tss 1\ot T^{\tss\ast}(v))
\,
R_{1,m+1}(u_1-v)\dots R_{m,m+1}(u_m-v)
\,
(\tss T_1(u_1)\dots T_m(u_m)\ot1\tss)\tss.
\end{gather*}
Here the subscript $m+1$ labels the last
tensor factor $\End\CC^N$, which comes after $\DY(\gl_N)$.
This relation can be proved by
using \eqref{3.99} repeatedly, i.e.\ $m$ times.
\epf

We have now established the following theorem explicitly decribing 
$\DY(\gl_N)\,$.

\bth\label{cor:doublerel}
The algebra $\DY(\gl_N)$ is generated by elements
$T_{ij}^{(r)},T_{ij}^{\tss(-r)}$ with $1\leqslant i,j\leqslant N$
and $r\geqslant 1$ subject only to the relations
\eqref{ternary},\eqref{3.13} and \eqref{3.111111}.
\eth

Note that
the relation \eqref{3.111111} is
equivalent to the collection of relations
\ben
(u-v)\ts [\tss T_{ij}(u),T^{\tss\ast}_{kl}(v)\tss]=
\sum_{m=1}^N
\Big(
\de_{jk}\,T_{im}(u)\,T^{\tss\ast}_{ml}(v)-
\de_{il}\,T^{\tss\ast}_{km}(v)\,T^{\tss\ast}_{mj}(u)
\Big)
\een
for all $i,j,k,l=1,\dots,N$.
We omit the proof of the equivalence, as
it is very similar to the proof of Proposition~\ref{prop:ternary}.
The last displayed relation can be rewritten as
\ben
[\tss T_{ij}(u),T^{\tss\ast}_{kl}(v)\tss]=
\sum_{p=0}^{\infty}\,
\sum_{m=1}^N\,
u^{-p-1}v^p\,
\Big(
\de_{jk}\,T_{im}(u)\,T^{\tss\ast}_{ml}(v)-
\de_{il}\,T^{\tss\ast}_{km}(v)\,T_{mj}(u)
\Big).
\een 
Expanding here the series in $u,v$ and
equating the coefficients at $u^{-r}v^{s-1}$ we get
\begin{gather*}
[\,T^{\tss(r)}_{ij}, T^{\tss(-s)}_{kl}]
=
\sum_{a=\max(1,r-s+1)}^{r}\,
\Big(
\de_{jk}
\Big(
\de_{a,r-s+1}\,T^{\tss(r-s)}_{il}+
\sum_{m=1}^N\,
T^{\tss(a-1)}_{im}\,T^{\tss(r-s-a)}_{ml}
\Big)
\\
\nonumber
-\,
\de_{il}
\Big(
\de_{a,r-s+1}\,T^{\tss(r-s)}_{kj}+
\sum_{m=1}^N\,
\,T^{\tss(r-s-a)}_{km}\,T^{\tss(a-1)}_{mj}
\Big)\Big)
\end{gather*}
for any indices $r,s\geqslant1$.
Here we keep to the notation $T_{ij}^{(0)}=\de_{ij}\,$.

We will complete this section with describing a bialgebra structure
on $\DY(\gl_N)$. The algebra $\DY(\gl_N)$ is generated by its two
subalgebras, $\Y(\gl_N)$ and $\Y^\ast(\gl_N)$.
We have already shown that the assignments 
\eqref{Delta} and \eqref{3.17} define comultiplications
on these two subalgebras,
while the assignments
$\ve:T(u)\mapsto1$ and $\ve:T^{\tss\ast}(v)\mapsto1$
define counit maps on them; see 
Propositions \ref{thm:hopf} and \ref{prop:dualhopf}.
Let us now replace the comultiplication $\Delta$ on $\Y^\ast(\gl_N)$
by its opposite comultiplication $\Delta'$.

\bpr
\label{prop:doublehopf}
The double Yangian $\DY(\gl_N)$ is a bialgebra
over\/ $\CC$ with the comultiplication defined by
extending $\Delta$ on $\Y(\gl_N)$ and $\Delta'$ on
$\Y^\ast(\gl_N)$, and with the counit defined by
mapping $T(u),T^{\tss\ast}(v)\mapsto1$.
\epr

\bpf
Using the equivalent form \eqref{3.111111}
of the defining relations \eqref{3.99},
the proof is similar to that
of the proof of Proposition~\ref{thm:hopf}. Here we omit the details.
\epf


\section{Filtration on the double Yangian}
\label{sec:dpbw}

In Section \ref{sec:filtr} we explained that
the associative algebra $\Y(\gl_N)$ can be regarded as a flat deformation
of the universal enveloping algebra $\U(\gl_N[z])$.
Our explanation was based on Proposition~\ref{prop:iso2}. 
In the present section we establish
an  analogue of that result for the double Yangian
$\DY(\gl_N)\,$. 

In order to do so, let us replace the descending
filtration on the algebra
$\Y^\ast(\gl_N)$ by an ascending filtration, such
that any generator $T_{ij}^{\ts(-r)}$ with $r\geqslant1$
has the degree $-\tss r\,$. Relative to this ascending filration
on $\Y^\ast(\gl_N)$, the subspace of elements of degree 
not more than $-\tss r$ coincides with the subspace of the elements
of degree not less than $r$ relative to the descending filtration. 
Let us now combine the ascending filtration on 
$\Y^\ast(\gl_N)$ with the ascending filtration on 
$\Y(\gl_N)$ used in Section \ref{sec:pbw}. That is, now
introduce an ascending $\ZZ\tss$-filtration on the algebra $\DY(\gl_N)$
by setting
\ben
\degpr T_{ij}^{\tss(r)}=r-1
\fand
\degpr T_{ij}^{\tss(-r)}=-\tss r
\een
for each index $r\geqslant1$.
Denote by $\grpr\DY(\gl_N)$ the corresponding $\ZZ\tss$-graded algebra.
Keeping to the notation of Section \ref{sec:pbw},
for any $r\geqslant1$
let $\Tt_{ij}^{\ts(r)}$ be the image of $T_{ij}^{\tss(r)}$
in the degree $r-1$ component of $\grpr\DY(\gl_N)\,$. 
Since we are now using an ascending filtration 
on $\Y^\ast(\gl_N)$ instead of the descending one,
for any $r\geqslant1$
we will denote by $\Tt_{ij}^{\ts(-r)}$ the image of $T_{ij}^{\tss(-r)}$
in the degree $-\tss r$ component of $\grpr\DY(\gl_N)$\,.
So $\Tt_{ij}^{\ts(-r)}$ now formally gets a new meaning,
which should not cause any confusion however.

\ble
\label{lem:gdrel}
In the graded algebra\/ $\grpr\DY(\gl_N)$
for any $r,s\geqslant1$ we have
\ben
[\,\Tt^{\ts(r)}_{ij},\Tt^{\ts(-s)}_{kl}]
\,=\,
\begin{cases}
\
\de_{kj}\,\Tt^{\ts(r-s)}_{il}
-\tss
\de_{il}\,\Tt^{\ts(r-s)}_{kj}
&\text{if}\quad r-s>0,
\\[8pt]
\
\de_{kj}\,\Tt^{\ts(r-s-1)}_{il}
-\tss
\de_{il}\,\Tt^{\ts(r-s-1)}_{kj}
&\text{if}\quad r-s\leqslant0.
\end{cases}
\een
\ele

\bpf
This follows from the relation displayed in
Section \ref{sec:doubleyang} last.
Indeed, relative to the ascending filtration on $\DY(\gl_N)$
the commutator at the left hand side of that relation
has the degree $r-s-1$ for any $r,s\geqslant1$. For $r-s>0$
the sum at the right hand side 
equals
\ben
\de_{jk}\,T^{\tss(r-s)}_{il}-\de_{il}\,T^{\tss(r-s)}_{kj}
\een
plus terms of degree not more that $r-s-2\,$. For $r-s=0$ 
that sum equals
\ben
\de_{jk}\,\big(\,\de_{il}+T^{\tss(-1)}_{il}\,\big)-
\de_{il}\,\big(\,\de_{kj}+T^{\tss(-1)}_{kj}\,\big)
=
\de_{jk}\,T^{\tss(-1)}_{il}-
\de_{il}\,T^{\tss(-1)}_{kj}
\een
plus terms of degree not more that $-2$. Finally, for $r-s<0$
that sum equals
\ben
\de_{jk}\,T^{\tss(r-s-1)}_{il}-\de_{il}\,T^{\tss(r-s-1)}_{kj}
\een
plus terms of degree not more that $r-s-2\,$.
\epf

The graded algebra $\grpr\DY(\gl_N)$
inherits from $\DY(\gl_N)$ a bialgebra structure, see
Proposition~\ref{prop:doublehopf}. Moreover
$\grpr\DY(\gl_N)$ is a Hopf algebra, see 
the remarks we made just before Proposition~\ref{prop:C4.3}.

\bpr
\label{thm:gdpbw}
The graded Hopf algebra $\grpr\DY(\gl_N)$ is isomorphic to
universal enveloping algebra $\U(\tss\gl_N[z,z^{-1}]\tss)\,$.
\epr

\bpf
Consider the subalgebras $\grpr\Y(\gl_N)$ and
$\grpr\Y^\ast(\gl_N)$ of the graded algebra $\grpr\DY(\gl_N)\,$.
We have an isomorphism \eqref{surhom} of graded algebras
defined by the assignments \eqref{asshom}.
Further, due to Lemma~\ref{lem:grel}
a surjective homomorphism
$$
\U(\tss z^{-1}\tss\gl_N[z^{-1}]\tss)\to\grpr\Y^\ast(\gl_N)
$$
can be defined by
$$
E_{ij}\,z^{-r}\mapsto\Tt^{\tss(-r)}_{ij}
\quad\text{for}\quad
r\geqslant1\,.
$$
Lemma~\ref{lem:gdrel} ensures that these two homomorphisms
extend to a homomorphism
\beql{dsurhom}
\U(\tss\gl_N[z,z^{-1}]\tss)\to\grpr\DY(\gl_N)\,.
\eeq
This homomorphism is surjective
and we will prove that it is injective as well.
Our proof will be similar to the proof 
of injectivity of the homomorphism
\eqref{surhom} given at the end of Section~\ref{sec:pbw}.
But now we will use Propositions~\ref{prop:C4.6} and \ref{prop:doublehopf}.

Take any finite linear combination $C$ of the products
\ben
(E_{\tss i_1j_1}\tss z^{\tss s_1})
\dots
(E_{\tss i_mj_m}\tss z^{\tss s_m})
\in\U(\tss\gl_N[z,z^{-1}]\tss)
\een
with certain coefficients
\ben
C_{\tss i_1j_1\dots i_mj_m}^{\,s_1\dots s_m}\in\CC
\een
where the indices $s_1\tss,\dots,s_m\in\ZZ$ and the number
$m\geqslant0$ may vary; the indices
$i_1,j_1,\dots,i_m,j_m\tss$ may vary as well.
Suppose $C\neq0$ as
an element of $\U(\gl_N[z,z^{-1}])\,$.
The algebra $\U(\gl_N[z,z^{-1}])$ comes with a natural $\ZZ$-grading 
such that
for any integer $s$ the generator $E_{ij}\tss z^s$ has the degree $s$.
The homomorphism \eqref{dsurhom} preserves this grading.
Without loss of generality, suppose that the element $C$ is
homogeneous of degree $d$ with respect to this grading. That is,
\ben
C_{\tss i_1j_1\dots i_mj_m}^{\,s_1\dots s_m}=0
\ \quad\text{if}\quad
s_1+\dots+s_m\neq d\tss.
\een
Now define the element $A\in\DY(\gl_N)$ as the sum
\ben
\sum_{s_1+\dots+s_m=d}
C_{\tss i_1j_1\dots i_mj_m}^{\,s_1\ldots s_m}\,
T_{i_1j_1}^{(r_1)}\dots T_{i_mj_m}^{(r_m)}
\een
where for every $k=1,\dots,m$ we set $r_k=s_k$ if $s_k<0\,$,
and $r_k=s_k+1$ if $s_k\geqslant0\,$.
Let $B$ be the image of $A$
in the $d\tss$-th component of the graded algebra $\grpr\DY(\gl_N)\,$.
The element $B$ coincides with the image of $C$
under the homomorhism \eqref{dsurhom}.

For any non-zero complex number $c$
the evaluation representation \eqref{bareval}
of the algebra $\U(\gl_N[z])$ can be extended
to a representation $\st_c$ of $\U(\gl_N[z,z^{-1}])$ so that 
\ben
\st_c:
E_{ij}\tss z^{s}\mapsto c^{\tss s}\tss e_{ij}
\,\quad\text{for any}\quad
s\in\ZZ\tss.
\een
Then by \eqref{3.666666} and \eqref{dualcovec} we have
\ben
\st_c(\tss E_{ij}\tss z^s\tss)=
\begin{cases}
\ \si_c(\tss T_{ij}^{(s)}\tss)
&\text{if}\quad s<0\,,
\\[2pt]
\ \si_c(\tss T_{ij}^{(s+1)}\tss)
&\text{if}\quad s\geqslant0\,.
\end{cases}
\een

Now let $c_1,\dots,c_n$ be any non-zero complex numbers.
Let $D\in(\End\CC^N)^{\ot n}$ be the image of $C$
under the tensor product of the representations
$\st_{c_1},\dots,\st_{c_n}$ of the algebra 
$\U(\gl_N[z,z^{-1}])\tss$. Denote by
$\si_{c_1\dots\tss c_n}$ the tensor product
of the representations $\si_{c_1},\dots,\si_{c_n}$
of the algebra $\DY(\gl_n)\tss$; here we use 
Proposition~\ref{prop:doublehopf}.
The image of $A\in\DY(\gl_N)$ under the representation
$\si_{c_1\dots\tss c_n}$ is a Laurent polynomial in $c_1,\dots,c_n\ts$.
The degree of this polynomial does not exceed $d\tss$, see
\eqref{Deltamatelem} and \eqref{Deltamadual}.
The sum of the terms of degree $d$
of this polynomial equals $D$, see the proof of Proposition~\ref{thm:P2.2}.

For any finite-dimensional Lie algebra
$\agot$ there is an analogue of Lemma~\ref{lem:L2.1} 
for $\agot\tss[z,z^{-1}]$ instead of $\agot\tss[z]$.
The proof of that analogue is similar to that of Lemma~\ref{lem:L2.1}
itself and is omitted here. Using that analogue,
we can choose $n$ and $c_1,\dots,c_n\neq0$ so that $D\neq0\tss$. 
Then $\degpr A=d$.
Indeed, if we had $\degpr A<d$ then the degree of
the Laurent polynomial $\si_{c_1\dots\tss c_n}(A)$ would be also less
then $d\tss$. This would contradict to the non-vanishing of 
$D$. By the definition of the element $B\in\grpr\Y(\gl_N)\tss$,
the equality $\degpr A=d$ means that $B\neq0$.
So the homomorphism \eqref{dsurhom} is injective.

Comparing the definitions \eqref{gradedhopf1},\eqref{gradedhopf2}
and \eqref{gradeddualhopf},\eqref{gradedduals}
with general definitions 
\eqref{standardhopf1},\eqref{standardhopf2} now completes
the proof of the proposition.
\epf

By applying the Poincar\'e--Birkhoff--Witt theorem 
\cite[Section~2.1]{d:ae}
to the current Lie algebra $\gl_N[z,z^{-1}]$ we now obtain its
analogue for the double Yangian $\DY(\gl_N)\,$.

\bth
\label{thm:dpbw}
Given any linear ordering of the set of
generators $T^{(r)}_{ij}$ and $T^{(-r)}_{ij}$
with $r\geqslant 1\,$,
any element of the algebra $\DY(\gl_N)$ can be uniquely
written as a linear combination of ordered monomials in these generators.
\eth

\bco
The defining homomorphisms of the algebras
$\Y(\gl_N)$ and\/ $\Y^\ast(\gl_N)$ to $\DY(\gl_N)$ are embeddings.
\eco

We will now use our ascending filtration on
$\DY(\gl_N)$ to show that in the initial definition of this algebra,
the relations \eqref{3.99} can be replaced by the relations
\beql{3.999}
\Delta(X)\,\tss\Rc
=
\Rc\,\tss\Delta'(X)
\quad\,\text{for every}\quad
X\in\Y(\gl_N)\,.
\eeq
Here $\Delta'$
is the comultiplication on $\Y(\gl_N)$
opposite to 
\eqref{Delta}. The infinite sums at both sides of the
relations \eqref{3.999} 
can be regarded as elements of the tensor product of
$\Y(\gl_N)$ and of the completion of $\DY(\gl_N)$ relative to
our ascending filtration. The completion of
$\Y^\ast(\gl_N)$ as a subalgebra of $\DY(\gl_N)$
then coincides with $\Y^\circ(\gl_N)\,$.

\bpr
Relations \eqref{3.999} in the algebra
$\DY(\gl_N)$ are equivalent to~\eqref{3.99}.
\epr

\bpf
Let $Y_1,Y_2,\dots$ be the basis of $\Y(\gl_N)$
from the proof of Proposition~\ref{prop:dualsys}.~Let
\ben
Y_p\,Y_q=\sum_{r=1}^\infty\, a_{pq}^{\tss r}\,Y_r
\quad\text{and}\quad
\De\tss(Y_r)\,=\sum_{p,q=1}^\infty b_{pq}^{\tss r}\,Y_p\ot Y_q
\een
so that $a_{pq}^{\tss r}\,,b_{pq}^{\tss r}\in\CC$ are 
the structure constants of the
bialgebra $\Y(\gl_N)$ relative to this basis.
Since the system of vectors 
$Y'_1,Y'_2,\dots$ of $\Y^\circ(\gl_N)$ is dual to the system
$Y_1,Y_2,\dots$ relative to the bialgebra pairing 
\eqref{pairing}, we also have the equalities
\ben
Y'_p\,Y'_q=\sum_{r=1}^\infty b_{pq}^{\tss r}\,Y_r'
\quad\text{and}\quad
\De\tss(Y'_r)=\sum_{p,q=1}^\infty a_{pq}^{\tss r}\,Y_p'\ot Y'_q\,.
\een
Here we extend the comultiplication 
$\Delta$ on $\Y^\ast(\gl_N)$ to
$\Y^\circ(\gl_N)$ as we did just after stating
Proposition \ref{prop:dualhopf}.

It suffices to take $X=Y_r$ with $r=1,2,\ldots$ in the relations 
\eqref{3.999}. Hence we get
$$
\sum_{p,q,s=1}^\infty b_{pq}^{\tss r}\,
(\tss Y_p\,Y_s')\ot(\tss Y_q\,Y_s)
\ =
\sum_{p,q,s=1}^\infty b_{pq}^{\tss r}\,
(\tss Y_s'\,Y_q)\ot(\tss Y_s\,Y_p)
$$
or
$$
\sum_{p,q,s,t=1}^\infty 
a_{qs}^{\tss t}\,b_{pq}^{\tss r}\,(\tss Y_p\,Y_s')\ot Y_t
\ =
\sum_{p,q,s,t=1}^\infty 
a_{sp}^{\tss t}\,b_{pq}^{\tss r}\,(\tss Y_s'\,Y_q)\ot Y_t
\,.
$$
So the relations \eqref{3.999} are equivalent to
the relations in our completion of $\DY(\gl_N)$
\beql{abpq}
\sum_{p,q,s=1}^\infty 
a_{qs}^{\tss t}\,b_{pq}^{\tss r}\,Y_p\,Y_s'
\ =
\sum_{p,q,s=1}^\infty 
a_{sp}^{\tss t}\,b_{pq}^{\tss r}\,Y_s'\,Y_q
\quad\text{where}\quad
r\tss,t=1,2,\ldots\,.
\eeq

The vectors $Y'_1,Y'_2,\dots$ have been determined by \eqref{gsum}
using a basis $Y^\ast_1,Y^\ast_2,\dots$ of $\Y^\ast(\gl_N)\,$. 
We also have the equalities
\beql{fsum}
Y^\ast_s=\sum_{r=1}^\infty\,f_{\tss rs}\,Y'_r
\eeq
where $f_{rs}=\langle\ts Y_r\tss,Y^\ast_s\ts\rangle\,$.
The matrix $[\,g_{rs}\tss]$ used in \eqref{gsum} is 
inverse to $[\,f_{rs}\tss]\,$.
Due to \eqref{gsum} and \eqref{fsum}
we can replace $W\in\Y^\ast(\gl_N)$ by $Y'_t\in\Y^\circ(\gl_N)$
with $t=1,2,\ldots\,$ in the relations \eqref{3.99}.
In this way we get
$$
\sum_{p,q,s=1}^\infty a_{pq}^{\tss t}\,
(\tss Y_s'\,Y'_p)\ot(\tss Y_s\,Y'_q)\ =
\sum_{p,q,s=1}^\infty a_{pq}^{\tss t}\,
(\tss Y'_q\,Y_s')\ot(\tss Y'_p\,Y_s)
$$
or
$$
\sum_{p,q,r,s=1}^\infty a_{pq}^{\tss t}\,b_{sp}^{\tss r}\,
Y'_r\ot(\tss Y_s\,Y'_q)\ =
\sum_{p,q,r,s=1}^\infty a_{pq}^{\tss t}\,\,b_{qs}^{\tss r}\,
Y'_r\ot(\tss Y'_p\,Y_s)\,.
$$
So the relations \eqref{3.99} are equivalent to
the relations in our completion of $\DY(\gl_N)$
\ben
\sum_{p,q,s=1}^\infty a_{pq}^{\tss t}\,b_{sp}^{\tss r}\,Y_s\,Y'_q\ =
\sum_{p,q,s=1}^\infty a_{pq}^{\tss t}\,\,b_{qs}^{\tss r}\,Y'_p\,Y_s
\quad\text{where}\quad
r\tss,t=1,2,\ldots\,.
\een
By cyclically permuting the summation indices in these relations we get
\eqref{abpq}.
\epf

\bco
\label{lastcor}
The coefficients of the series $Z(u)$ lie in the centre of\/ 
$\DY(\gl_N)\,$.
\eco

\bpf
The coefficients of 
$Z(u)$ lie in the centre of 
$\Y(\gl_N)$ by Lemma \ref{zu}. To prove that they commute with 
the elements of $\Y^\ast(\gl_N)$ as a subalgebra of $\DY(\gl_N)$
let us substitute the series $Z(u)$ for $X\in\Y(\gl_N)$ in \eqref{3.999}.
Due to Proposition \ref{lz} we get
$$
\sum_{s=1}^\infty\,(\tss Z(u)\,Y'_s\tss)\ot (\tss Z(u)\,Y_s)
\,=\,
\sum_{s=1}^\infty\,(\tss Y'_s\,Z(u))\ot (\tss Y_s\,Z(u))\,.
$$
As the coefficients of $Z(u)$ are central in $\Y(\gl_N)\,$,
dividing 
this 
by $1\ot Z(u)$ yields
$$
\sum_{s=1}^\infty\,(\tss Z(u)\,Y'_s\tss)\ot Y_s
\,=\,
\sum_{s=1}^\infty\,(\tss Y'_s\,Z(u))\ot Y_s\,.
$$
It follows that the coefficients of $Z(u)$ commute with every $Y'_s$ in
our completion of the algebra $\DY(\gl_N)\,$.
By using the relations \eqref{fsum} we now get the corollary.
\epf

Now consider the series $Z^{\tss\circ}(v)$ appearing in
Lemma \ref{zc}. Arguing as in the proof of the Corollary \ref{lastcor},
but using the relations \eqref{3.99} and Proposition \ref{zzcc}
instead of the relations \eqref{3.999} and Proposition \ref{lz},
we can show that the coefficients of 
$Z^{\tss\circ}(v)$ belong to the centre of our completion of the algebra  
$\DY(\gl_N)\,$.
However, in general these coefficients do not belong to the algebra
$\DY(\gl_N)$ itself, see Section \ref{sec:dualyang} again. 

Our completion of the algebra $\DY(\gl_N)$ can also be used 
to rewrite the relations \eqref{3.13} and \eqref{3.111111} 
similarly to \eqref{ternary}. Take the element
$T^{\tss\natural}(v)$ inverse to $T^{\tss\ast}(v)\,$. 
In the notation analogous to \eqref{tbast} 
the equality \eqref{3.13} of series in $u$ and $v$
with coefficients in $\Y^\ast(\gl_N)\ot(\End\CC^N)^{\ot 2}$
can be then rewritten as the equality
$$
R(u-v)\ts T^{\tss\natural}_1(u)\ts T^{\tss\natural}_2(v)=
T^{\tss\natural}_2(v)\ts T^{\tss\natural}_1(u)\ts R(u-v)
$$
of series 
with coefficients in $\Y^\circ(\gl_N)\ot(\End\CC^N)^{\ot 2}\,$.
The \eqref{3.111111} can be rewritten~as
$$
R(u-v)\,(\tss T(u)\ot1\tss)\,(\tss1\ot T^{\tss\natural}(v))
=
(\tss1\ot T^{\tss\natural}(v))\,(\tss T(u)\ot1\tss)\,R(u-v)\,.
$$



\end{document}